\documentclass[sn-mathphys-num]{sn-jnl}

\usepackage{graphicx}%
\usepackage{multirow}%
\usepackage{amsmath,amssymb,amsfonts}%
\usepackage{amsthm}%
\usepackage{mathrsfs}%
\usepackage[title]{appendix}%
\usepackage{xcolor}%
\usepackage{textcomp}%
\usepackage{manyfoot}%
\usepackage{booktabs}%
\usepackage{algorithm}%
\usepackage{algorithmicx}%
\usepackage{algpseudocode}%
\usepackage{listings}%

\usepackage{mathtools}
\usepackage{physics}
\usepackage{hyperref}

\usepackage{enumitem}
\makeatletter
\def\dicho#1{\expandafter\@dicho\csname c@#1\endcsname}
\def\@dicho#1{\ifnum#1>1 or \else\fi(\@Roman#1)}
\makeatother
\AddEnumerateCounter{\dicho}{\@dicho}{or (III)}
\newlist{dichotomy}{enumerate}{1}
\setlist[dichotomy]{label=\dicho*,leftmargin=1.5cm}

\DeclarePairedDelimiter\floor{\lfloor}{\rfloor}

\newcommand{\tsum}{\mathop{\textstyle\sum}\nolimits}
\newcommand{\Z}{\mathbb{Z}}
\newcommand{\ab}{\underline{a}}

\newcommand{\smvee}{\raise0.4ex\hbox{$\scriptscriptstyle\vee$}}

\newcommand{\EO}{Ekedahl-Oort }

\newcommand{\ux}{\underline{x}}
\newcommand{\F}{\mathbb{F}}

\newcommand{\OO}{\mathcal{O}}

\newcommand{\M}{\mathcal{M}}

\newcommand{\Q}{\mathbb{Q}}
\newcommand{\sh}{\text{Sh}}

\newcommand{\Sh}{{\rm Sh}}
\newcommand{\Spec}{\text{Spec}}

\newcommand{\im}{\text{Im}}

\newcommand{\ptnk}{\lfloor{p\langle\frac{\tau a(k)}{m}\rangle}\rfloor}

\newcommand{\C}{\mathbb{C}}
\newcommand{\fpb}{\overline{\mathbb{F}}_p}
\newcommand{\sigp}{\sigma_p}

\newcommand{\Po}{\mathbb{P}^1}

\newcommand{\gra}{(G,r,\underline{a})}

\newcommand{\shgf}{\text{Sh}(\D)}

\newcommand{\ts}{\tau^*}

\newcommand{\tp}{\tau'}

\newcommand{\tps}{\tau'^*}
\newcommand{\ti}{\tau_i}

\newcommand{\phit}{\phi_{\tau}}
\newcommand{\psit}{\psi_{\tau}}

\newcommand{\tg}{\mathcal{T}_{G}}

\newcommand{\pbtar}{\lfloor{p\langle\tfrac{\tau a(r)}{m}\rangle}\rfloor}
\newcommand{\pbtr}{\lfloor{p\langle\frac{p^i \tau a(r)}{m}\rangle}\rfloor}
\newcommand{\pbtrm}{\lfloor{p\langle\frac{p^i \tau a(r-1)}{m}\rangle}\rfloor}

\newcommand{\jj}{\mathfrak{J}}

\newcommand{\hy}{H^1(C,\OO_C)}
\newcommand{\hyd}{H^0(C,\Omega_C)}

\newcommand{\pbtt}{\lfloor{p\langle\tfrac{\tau a(2)}{m}\rangle}\rfloor}
\newcommand{\pbto}{\lfloor{p\langle\tfrac{\tau a(1)}{m}\rangle}\rfloor}
\newcommand{\pbtf}{\lfloor{p\langle\tfrac{\tau a(4)}{m}\rangle}\rfloor}
\newcommand{\pbtd}{\lfloor{p\langle\tfrac{\tau a(3)}{m}\rangle}\rfloor}

\newcommand{\pbttinl}{\lfloor{p\langle\frac{\tau a(2)}{m}\rangle}\rfloor}
\newcommand{\pbtoinl}{\lfloor{p\langle\frac{\tau a(1)}{m}\rangle}\rfloor}
\newcommand{\pbtfinl}{\lfloor{p\langle\frac{\tau a(4)}{m}\rangle}\rfloor}
\newcommand{\pbtdinl}{\lfloor{p\langle\frac{\tau a(3)}{m}\rangle}\rfloor}

\newcommand{\pby}{\lfloor{p\langle\frac{ia(1)}{m}\rangle}\rfloor}

\newcommand{\ptk}{\lfloor{p\langle\tfrac{\tau a(k)}{m}\rangle}\rfloor}

\newcommand{\pbr}{\lfloor{p\langle\frac{ia(r)}{m}\rangle}\rfloor}

\newcommand{\pbtit}{\lfloor{p\langle\tfrac{\tau_i a(2)}{m}\rangle}\rfloor}
\newcommand{\pbtid}{\lfloor{p\langle\tfrac{\tau_i a(3)}{m}\rangle}\rfloor}
\newcommand{\pbtio}{\lfloor{p\langle\tfrac{\tau_i a(1)}{m}\rangle}\rfloor}
\newcommand{\pbtif}{\lfloor{p\langle\tfrac{\tau_i a(4)}{m}\rangle}\rfloor}

\newcommand{\pbk}{\floor*
{p\langle\frac{ia_k}{m}\rangle}}

\newcommand{\fts}{f(\ts)}

\newcommand{\sti}{s_{\tau_i}}

\newcommand{\au}{\underline{a}}

\newcommand{\gal}{\text{Gal}}

\newcommand{\tns}{\tau_0^{*}}
\newcommand{\fb}{\underline{f}}

\newcommand{\Hom}{\text{Hom}}

\newcommand{\orho}{\OO_{\tau}}

\newcommand{\mlra}{\mathcal{M}(m,r,\au)}

\newcommand{\mgra}{\mathcal{M}(G,r,\au)}

\newcommand{\A}{\mathcal{A}}

\newcommand{\ub}{\underline{b}}

\newcommand{\fp}{\mathbb{F}_p}
\newcommand{\mgrab}{\mgra_{\fpb}}

\newcommand{\tn}{\tau_0}

\newcommand{\st}{s_\tau}

\newcommand{\psitjjp}{\psi_\tau(j',j)}

\newcommand{\D}{\mathcal{D}}

\newcommand{\bu}{(b_1, b_2, \dots, b_r)}
\newcommand{\bun}{(b_1, b_2, \dots, b_r ; N)}
\newcommand{\bunj}{(b_1, b_2, \dots, b_r ; N;j')}

\newcommand{\bifor}{\frac{1+\lfloor p \langle\frac{\tau_i a(4)}{m}\rangle\rfloor}{1+\lfloor p\langle\frac{\tau_i a(1)}{m}\rangle\rfloor}}

\theoremstyle{thmstyleone}%
\newtheorem{theorem}{Theorem}
\newtheorem{proposition}[theorem]{Proposition}%

\theoremstyle{thmstyletwo}%
\newtheorem{remark}{Remark}%
\newtheorem{corollary}{Corollary}
\newtheorem{lemma}{Lemma}

\theoremstyle{thmstylethree}%
\newtheorem{definition}{Definition}%
\newtheorem{notation}{Notation}

\begin{document}

\title[Non-$p$-ordinary smooth abelian covers of $\Po$]{Non-$p$-ordinary smooth abelian covers of $\Po$}


\author*[1]{\fnm{Yuxin} \sur{Lin}\orcid{https://orcid.org/0000-0001-9230-7728}}\email{yuxinlin@caltech.edu}

\author[1]{\fnm{Elena} \sur{Mantovan}}\email{mantovan@caltech.edu}
\equalcont{These authors contributed equally to this work.}

\author[2]{\fnm{Deepesh} \sur{Singhal}\orcid{https://orcid.org/0000-0003-1317-5523}}\email{singhald@uci.edu}
\equalcont{These authors contributed equally to this work.}

\affil*[1]{\orgdiv{Department of Mathematics}, \orgname{California Institute of Technology}, \orgaddress{\street{1200 E. California Blvd.}, \city{Pasadena}, \postcode{91125}, \state{CA}, \country{USA}}}

\affil[2]{\orgdiv{Department of Mathematics}, \orgname{University of California, Irvine}, \orgaddress{\street{120 Aldrich Hall}, \city{Irvine}, \postcode{92697}, \state{CA}, \country{USA}}}



\abstract{Given a family of odd abelian covers of $\Po$ 
and a prime $p$ of good reduction, in \cite{lin2024abelian}, under some assumptions, we computed the generic Newton polygon (resp. Ekedahl-Oort type) in the family (called {\em $p$-ordinary}). In this paper, we investigate the existence of non-$p$-ordinary smooth curves in the family.
In particular, under some restrictions, we show that when $p$ is sufficiently large, the complement of the $p$-ordinary locus is always nonempty. For $1$-dimensional families satisfying additional auxiliary conditions, we obtain a lower bound for the number of non-$p$-ordinary  smooth curves.  In specific instances,  the above general statement can be improved; as an example, for families of covers of degree at most 7,  we establish the non-emptiness of certain non-$p$-ordinary Newton/Ekedahl-Oort strata  (called {\em almost $p$-ordinary}). Our method relies on further study of the extended Hasse-Witt matrix introduced by Moonen in \cite{moonen2022computing}, and  initiated in \cite{lin2024abelian}, and on known results about the geometry of the mod-$p$ reduction of Shimura varieties of PEL type.   }




\maketitle

\section{Introduction}
This paper is motivated by the arithmetic Schottky problem in positive characteristics, which investigates the mod-$p$ invariants of abelian varieties that occur for Jacobians of smooth curves. (We refer to \cite{PriesSurvey} and \cite{katzExpository} for an introduction to this problems, and related open questions).
We focus our attention to the case of Jacobians of abelian covers of the projective line $\Po$.

Let $G$ be a finite abelian group, of odd size $|G|$ and exponent $e$. We denote by $E=\Q [e^{\frac{2\pi i}{e}}]\subset \C$ the $e$-th cyclotomic field, and by $\OO_E$ its ring of integers.   The Hurwitz space $\mathcal{M}(G)$ of $G$-covers of $\Po$ (as defined in \cite[Section 2.1-2.2]{achter2007integral}) is a smooth and proper Deligne-Mumford stack defined over $\Spec(\OO_E[\frac{1}{e}])$, that parameterizes $G$-covers of $\Po$. On each irreducible component of $\mathcal{M}(G)$, the monodromy datum of the parameterized covers is constant, and conversely each monodromy datum uniquely determines an irreducible component.  In the following, we denote by $\mgra$, the irreducible component of the Hurwitz space associated with a monodromy datum $(G,r,\au)$, and by $g$ be the genus of the curves it parameterized.\footnote{In the notation $(G,r, \ab)$, $r$ denotes the number of branched points and $\ab=(a(1), a(2), \dots, a(r))$ is a $r$-tuple of elements in $G$, where for each $1\leq i\leq r$ the element $a(i) \in G$ is the local monodromy at the $i$-th branched point. } 
Let $\mathcal{A}_g$ denote the moduli space of principally polarized abelian varieties of dimension g. Over $\C$, the image 
of  $\mgra_{\C}$  under the Torelli  map is naturally contained in a special subvariety $S(G,r,\ab)$  of $\mathcal{A}_{g, \C}$. 
By definition (\cite{deligne1986monodromy}), 
$S(G,r,\au)$ is an irreducible component of a PEL-type moduli space $\Sh(\mathcal{D})_{\C}$, where $\Sh(\D)$ is a smooth proper Deligne-Mumford stack defined over $\Spec(\OO_E[\frac{1}{e}])$; furthermore, the Torelli morphism $\mgra_\C\to S(\gra)$ extends to a morphism $\mgra\to \Sh(\D)$ defined over $\Spec(\OO_E[\frac{1}{e}])$ (here, $\D=\D(G,r,\au)$ denotes the integral PEL-type moduli datum associated with the monodromy datum $(G,r,\au)$).

 For a prime $\mathfrak{p}$ of $\OO_E$ above a rational prime $p$, $p\nmid e$, we denote by  $\overline{\F}_p$ an algebraic closure of its residue field. The induced map from $\mgra_{\overline{\F}_p}$ to $\Sh(\D)_{\overline{\F}_p}$ gives rise to explicit constraints on the Ekedahl-Oort types and the Newton polygons of the curves parameterized by $\mgra_{\overline{\F}_p}$.
For example, the lowest Ekedahl–Oort type and Newton polygon occurring on $\Sh(\D)_{\overline{\F}_p}$ are natural lower bounds for those occurring on $\M(G,r,\au)_{\overline{\F}_p}$. By \cite[Theorem 1.3.7]{moonen2004serre},  the unique open Ekedahl–Oort and Newton strata of $\Sh(\D)_{\overline{\F}_p}$ agree. We refer to the pullback to $\M(G,r,\au)_{\overline{\F}_p}$ of this unique open and dense stratum of $\Sh(\D)_{\overline{\F}_p} $ as the $\D$-ordinary locus (at $\mathfrak{p}$).

In \cite[Theorem 1.1] {lin2024abelian}, for $p> (r-2)|G|$ , we showed if $r\leq 5$ then the $\D$-ordinary locus of $\M(G,r,\au)_{\overline{\F}_p}$ is non-empty. In this paper, we investigate whether its complement, the non-$\D$-ordinary locus of $\M(G,r,\au)_{\overline{\F}_p}$, is also non-empty; we give a positive answer in the following setting. 
%
%
More precisely, we prove the following result. 

\begin{theorem}\label{main theorem}
Let $(G,r,\au)$ be a monodromy datum, with $G$ a finite abelian group of odd size. Denote $\tg$ the character group of $G$, and  $f:\tg\to \Z_{\geq 0}$ the signature type of $(G,r,\au)$ (Definition \ref{signature type of monodromy}). 

Assume (A1): there exists $\tau \in \tg$ such that $f(\tau)=1$ and $g(\tau) \geq 2$. 

Then,  the non-$\D$-ordinary locus of $\mgra_{\overline{\F}_p}$ is non-empty,  for any prime $p > (3\binom{r}{2}+r-2)|G|$, and for  $p>3|G|$ if $r=4$.
\end{theorem}

The following corollary identifies monodromy data satisfying Assumption (A1) in Theorem \ref{main theorem}.
\begin{corollary}\label{cor: cases of main theorem}
Let  $(G,r,\au)$ be an odd abelian monodromy datum.  Assume one of the following holds:
\begin{enumerate}
\item $4\leq r\leq 5$.
\item for a proper subgroup $H\leq G$: $4\leq |\{1\leq i\leq r: a(i)\not \in H\}|\leq 5$. 
\item $G=\Z/m\Z$, and $\sum_i a(i)\in \{m, 2m, (r-2)m, (r-1)m\}$.
\end{enumerate}

Then, assumption (A1) of Theorem \ref{main theorem} is satisfied; hence 
 the non-$\D$-ordinary locus of $\mgra_{\overline{\F}_p}$ is non-empty, for any prime $p>(3\binom{r}{2}+r-2)|G|$, and for $p>3|G|$ if $r=4$.
\end{corollary}

When $r=4$,  the non-$\D$-ordinary locus has dimension 0. Under a stronger assumption,  we obtain a lower bound on the number of non-$\D$-ordinary points. (Note that when $r=4$ then $g(\tau)\geq 2$ implies $g(\tau)=2$; see Lemma \ref{lem: ft+fts<=r-2}.)


\begin{theorem}\label{prop: main theorem for r eq 4}
Let $(G, 4, \au)$ be an odd abelian monodromy datum, and $p$ a prime, $p\nmid |G|$.

Assume (A2): there exists $\tau \in \tg$ such that $f(\tau)=1$, $g(\tau)=2$, and $f(p\tau)=1$.

Then, for $p\geq 3|G|$,  the number of non-$\D$-ordinary points in  $\M(m, 4,\au)(\fpb)$ is at least $\lfloor{\frac{p}{|G|}-1}\rfloor$. 
\end{theorem}

In the case when $G$ is cyclic, and the $\D$-ordinary Newton polygon is ordinary (e.g., if $p\equiv 1 \bmod |G|$), a more precise formula for the number of non-ordinary curves in the family is given by Cavalieri-Pries in \cite[Main Theorem]{cavalieri2023mass}. Also, we have recently learned that, in the case when $m$ is prime and the $p$-rank of the $\D$-ordinary Newton polygon is positive, work in progress of Mohao Yi proves the existence and bounds the number of points of non-maximal $p$-rank.

As the non-$\D$-ordinary locus of $\mgra_{\overline{\F}_p}$ is naturally stratified into smaller locally closed Ekedahl-Oort and Newton strata, it is naturally to ask whether Theorems \ref{main theorem} and \ref{prop: main theorem for r eq 4} can be improved to identify the non-empty (non-$\D$-ordinary) strata in $\mgra_{\overline{\F}_p}$.
This can be achieved in specific instances.
As an example, we establish the following result for the cyclic monodromy datum $(7,4,(1,1,2,3))$. 

\begin{proposition}\label{prop: special case}
\hfill
 Let $p$ be a prime,  $p\neq 7$. 
Denote $\M=\M(\Z/7\Z,4,(3,2,1,1))$, and $\Sh=\Sh(\D (\Z/7\Z,4,(3,2,1,1))$. 

If $p\not \equiv 1,6\pmod{7}$:  the Newton and Ekedahl-Oort stratifications of $\Sh_{\fpb}$ agree, there are two strata, and both have non-empty intersection with the image of $\M_{\overline{\F}_p}$ under the Torelli map.

If $p\equiv 1,6\pmod {7}$:  the Ekedahl-Oort stratification of $\Sh_{\fpb}$ is a refinement of the Newton stratification;  there are three Newton strata, and four Ekedahl-Oort strata (the Newton stratum of dimension 1 is the union of the two Ekedahl-Oort strata of dimension 1);  all strata have non-empty intersection with the image of $\M_{\overline{\F}_p}$ except:
\begin{enumerate}

\item 
both \EO strata of dimension 1 (and hence also the Newton stratum of dimension 1) for $p=13$;
\item  
the closed stratum of dimension 0 for some but not all  $p\equiv 1\pmod{7}$. 

\end{enumerate}

\end{proposition}

\begin{remark}
When the dimension of $\M_g$ is smaller than the codimension of a 
Newton stratum $\A_g[\nu]$ in $\A_g$, their (non-empty) intersection is considered unlikely (see \cite[Section 8]{li2019newton}). For instance, if $g>8$ then dimension of $\M_g$  is smaller than the codimension of the supersingular locus of $\A_g$. 

By \cite[Theorem 4.1]{de2000purity}, if $g>8$, then either \begin{itemize}
\item (O1) there exists a Newton polygon $\nu$ such that the corresponding (locally closed) Newton stratum $\M_g[\nu]$ is empty, or
\item (O2) there exist two Newton polygon $\nu'>\nu$ (here, $>$ means "lies above") and a connected component $Z$ of $\M_g[\nu]$ such that the intersection of the closure of $Z$ with $\M_g[\nu']$ is empty.
\footnote{For comparison, 
by \cite[Theorem 6.2, Remark 6.6]{Oort2000}, for any
two Newton polygon $\nu'>\nu$ and any connected component $W$ of $\A_g[\nu]$, the intersection of the closure of $W$ with $\A_g[\nu']$  
is always non-empty.}
\end{itemize}

Within the context of this paper, that is, for the variant of these
questions for families of abelian cover of $\Po$ and the associated PEL-type moduli spaces, Proposition 
\ref{prop: special case} include instances of both (O1), for $p=13$ and some $p\equiv 1\pmod{7}$  (see Remark \ref{some p}), and (O2), for $p\equiv 6\pmod{7}$.
\footnote{
 While instances of unlikely intersections  were previously established (see \cite[Section 8]{li2019newton} and  \cite[Remark 8.6]{lin2024abelian}), the latter two phenomena were not.}

\end{remark}

More generally, in Proposition \ref{prop: other fam}, we prove that, given a cyclic monodromy data $(m,4,\au)$ and a prime $p$, $p>3m$, if $m=5$ or $m=7$ then all Ekedahl-Oort strata of $\Sh(\D)_{\overline{\F}_p}$ of codimension 1 (these are called almost-$p$-ordinary), and  the unique Newton stratum  of $\Sh(\D)_{\overline{\F}_p}$ of codimension 1, have non-empty intersection with the image of  $\mgra_{\overline{\F}_p}$ under the Torelli map, unless the monodromy datum is equivalent to $(7,4, (1,6,c,7-c))$, for some $1\leq c\leq 6$.

\begin{remark}
    As in \cite[Remark 8.3]{lin2024abelian}, by  the results in \cite[Remark 8.5]{li2019newton}, \cite[Theorem 1.1]{lin2024abelian} and Theorem \ref{main theorem} combined give rise to infinitely many (inductive systems of) cyclic monodromy data, where degree $m$ is constant and the number of branched points grows unbounded, such that for $p>3m$  the corresponding families of covers of $\Po$ contains both $\D$-ordinary and non-$\D$-ordinary (resp. almost-$\D$-ordinary if  $m=5$ or $m=7$) smooth curves over $\fpb$. 
\end{remark}

The results in this paper rely on the study of the generalized Hasse--Witt invariant on $\mgrab$ constructed in \cite{lin2024abelian}.\footnote{In \cite{moonen2022computing}, Moonen introduced the notion of Hasse-Witt triple for a smooth curve over $\fpb$, which recovers the Ekedahl-Oort type of the Jacobian of the curve, and generalizes the classical Hasse-Witt invariant.   
In \cite{lin2024abelian}, we computed the Hasse-Witt triple  of odd abelian covers of the projective line, and give a criterium for $\D$-ordinariness in terms of the entries of Hasse-Witt triples, as the non-vanishing of a certain characteristic polyonomial, the generalized Hasse-Witt invariant.}  
Its vanishing locus corresponds to (possibly singular) non-$\D$-ordinary curves in the family; hence, we establish the existence of a smooth non-$\D$-ordinary curve by showing that the non-$\D$-ordinary locus is not fully contained in the boundary of the family (the locus parametrizing non-smooth covers). 
 When $r=4$, the generalized Hasse--Witt invariant can be regarded as a single variable polynomial with coefficients in $\fp$.  In the classical case of modular curve, it is a separable polynomial whose degree grows with $p$, and  its the vanishing locus hence does not consists solely of cusps. In our setting, the generalized Hasse--Witt invariant factors as a product of partial invariants, their degrees depending both on $p$ and the monodromy datum. We prove Theorem \ref{prop: main theorem for r eq 4} by identifying a non-trivial separable factor.\footnote{In general, the generalized Hasse--Witt invariant is a sum of the products of polynomials whose coefficients and degrees depend on the monodromy datum. While we are able to show that each of these summand is separable, their sum is generally not.} 
 By construction, each  almost-$\D$-ordinary \EO stratum is the vanishing locus of a partial Hasse--Witt invariant. Hence, in specific instances, by verifying the existence of non-common roots, we establish the non-emptiness of the corresponding almost-$\D$-ordinary \EO strata, and subsequently also of codimension 1 Newton strata (when it fully contains the almost-$\D$-ordinary \EO strata).

The paper is organized as follows. 
In Section \ref{pre}, we recall the notions of Ekedahl--Oort type, Newton polygon,  $\D$-ordinariness, Hasse--Witt triple and reduce to the study of cyclic covers. In Section \ref{sec: hasse invariant one dim}, we compute $\D$-ordinary Hasse--Witt invariant for families with four branched points, and compute its order of vanishing at $0$ and $1$ (that is, at the points representing non-smooth curves in the families). We prove Theorems \ref{main theorem} and \ref{prop: main theorem for r eq 4}  first for covers branched at four points  in Section \ref{sec:fourpoints},  and then for covers branched at five or more points in Section \ref{sec: main case}, where we also prove Corollary \ref{cor: cases of main theorem}. In Section \ref{sec: example of a family}, we prove Proposition \ref{prop: special case} and Proposition \ref{prop: other fam} for families  of cyclic covers of $\Po$, branched at four points, and of degree equal to either 5 or 7. In Appendix \ref{sec: separability}, we establish the properties of the polynomials $f(a,b,c)$ used in Sections \ref{sec: hasse invariant one dim} and \ref{sec: example of a family}. Appendix \ref{app: technical lemmas} contains the technical lemmas on the class of polynomials appearing in Section \ref{sec: main case}.

\section{Preliminaries}\label{pre}

\subsection{Notations}\label{tg}
Let $G$ be an abelian group, of size $|G|$. We assume throughout this paper that $|G|$ is odd. Denote $\tg={\rm Hom}(G,\C^\times)$ the  character group of $G$. For each $\tau\in\tg$, we definite its complex conjugate $\tau^*\in \tg$, by $\tau^*(g)=\tau(g^{-1})$ for all $g\in G$.

Let $p$ be a rational prime. We assume $p\nmid |G|$. We fix an algebraic closure $\overline{\Q}_p$ of $\Q_p$, together with  an isomorphism $\iota:\C\simeq {\overline{\Q}}_p$.
We denote $\overline{\Q}^{\rm un}_p$ the maximal unramified subfield of $\overline{\Q}_p$, $\fpb$ its residue field, and $\sigp\in \gal(\overline{\Q}^{\rm un}_p/\Q_p)$ the lift of the Frobenius element in $\gal(\fpb/\mathbb{F}_p)$.  

Let $e={\rm exp} (G)$ denote the exponent of $G$. Set  $E=\Q[\zeta_e]\subseteq \C$  where $\zeta_e=e^{2\pi i/e}\in \C$, and write $\OO_E$ for the ring of integer of $E$.   By assumption, $p\nmid e$; 
via $\iota$, we realize $E\subseteq \overline{\Q}_p^{\rm un}$; this yields a prime $\mathfrak{p}$ of $E$ above $p$ and an inclusion of the residue field of $\OO_E$ at $\mathfrak{p}$ in  $\fpb$.

Via $\iota$, we identify $\tg\simeq \Hom(G,\overline{\Q}^{\rm un, \cross}_{p})$, 
 and define an action of $\sigp$ on $\tg$ by composition. That is, $\tau^{\sigp}(x)=\sigp(\tau(x))=\tau(x)^p$.
 For each $\tau \in \tg$, we denote by $\orho$  the orbit of $\tau$ under Frobenius, and  by  $\mathfrak{O}_G$ the set of Frobenius orbits of $\tg$. 
For simplicity, given $\tau\in\tg$, we denote $\tau^{\sigp}$ by $p\tau$ and its orbit $\OO_{\tau}=\{\tau,p\tau,\dots,p^{l(\OO_{\tau})-1}\tau\}$, where $l(\OO_\tau)=|\OO_\tau|$.  Since $p(p^{l(\OO_\tau)-1}\tau)=\tau$, we also write $p^{l(\OO_\tau)-1}\tau$ as $\frac{\tau}{p}$.   When there is no ambiguity, we  write $\OO=\OO_\tau$ and $l=l(\OO)$.

When $G$ is a cyclic group of size $m$ for some odd $m$, we identify $\tg=\Z/m\Z$ as follows: for $i \in \Z/m\Z$, denote by $\tau_i$ the character  of $G$ such that $\tau_i(1)=\zeta_m^i$.

\subsection{The Hurwitz space of abelian covers of $\Po$} 

We briefly recall the definition of the Hurwitz space of $G$-covers of the projective line $\Po$, for $G$ a finite abelian group. We follow \cite[Section 2.2]{lin2024abelian}.

\begin{definition}\label{monodromy datum} 
Let $G$ be a finite abelian group. Consider a triple $(G,r,\au)$, where $r \geq 3$ is an integer and $\au=(a(1), \dots, a(r))$ is a $r$-tuple in $G^r$. The triple $(G,r,\au)$ is a \textbf{monodromy datum} if it satisfies the following conditions:
\begin{enumerate}
    \item $a(i) \neq 0$ in $G$,
    \item $a(1), \dots, a(r)$ generate $G$,
    \item $\sum_{i=1}^r a(i) = 0$ in $G$.
\end{enumerate}
We say that two monodromy data $(G,r,\au)$ and $(G,s,\ub)$ are {\bf equivalent} if $r=s$ and $\au, \ub \in G^r$ are in the same orbit under ${\rm Aut}(G) \times {\rm Sym}_r$. That is, if there exists $(\tau,\sigma)\in {\rm Aut}(G) \times {\rm Sym}_r$  such that $b(i)=\tau(a(\sigma(i)))$, for all $1\leq i\leq r$.
   
\end{definition}
For $C$ a smooth projective curve defined over an algebraically closed field $k$ of characteristic $p$, 
if $C$ is a $G$-cover of $\Po$ and $p\nmid |G|$, then the monodromy datum of $C$ is $(G,r,\au)$ where  $r$ is the number of branched points, and $\au$ is the  inertia type of $C$. We refer to \cite[Section 2.2]{lin2024abelian} for the definition of $\au$.

Given a monodromy datum $(G,r,\au)$, we denote by $\mgra$ the moduli stack of $G$-covers of $\Po$, with this monodromy datum.
It is a smooth Deligne--Mumford stack defined over $\mathbb Z[\frac{1}{|G|}]$.

\begin{definition}\label{signature type of monodromy}
    Let $C$ be $G$-cover of $\Po$ defined over $\C$, with monodromy datum $(G,r,\au)$. Set ${\rm V}={\rm H}^1(C(\C),\Q)$, and denote the Hodge decomposition of ${\rm V}$ by ${\rm V}\otimes_\Q\C={\rm V}^+\oplus {\rm V}^-$, where ${\rm V}^+={\rm H}^0(C,\Omega^1)$ and ${\rm V}^-={\rm H}^1(C,\OO)$. 
    
    Then, $C$ has genus $g=\dim V^+$, and its \textbf{signature type} $\fb$ is defined as
     $f :\tg\to \Z $, $ f(\tau)=\dim {{\rm V}}^{+}_{\tau},$ where for each character $\tau\in \tg$, we denote by ${\rm V}_\tau^+$ the subspace of ${\rm V}^+$ of weight $\tau$. 
    In particular,  $g=\sum_{\tau \in \tg}f(\tau)$.

\end{definition}

 In Section \ref{sec: reduction}, 
we reduce Theorems \ref{main theorem} and \ref{prop: main theorem for r eq 4} from the case of abelian covers to that of cyclic covers. Hence, in the following,  we focus on the case of $G$  cyclic.

Assume $G$ is cyclic, of size $m$ (odd). In this case, we also write the monodromy datum $(G, r,\au)$ as $(m, r, \au)$, where for $1\leq i\leq r$, $a(i)\in\Z$, $1\leq a(i)\leq m-1$.
Then, the curves parametrized by $\mlra$ are given by the normalization of the equation
\[y^m=\prod_{k=1}^r(x-x_k)^{a(k)},\]
where $x_1, \dots, x_r$ are pairwise distinct.
Let $C_{\ux}$ denote the curve associated with the point $\ux=(x_1,\dots,x_r) $ of $ \mlra$.
The genus $g$ and the signature $\fb$ are independent of the choice of the point $\ux\in \M(m,r,\au) (\C)$, and can be computed explicitly from the monodromy datum via the Hurwitz-Chevalley-Weil formula (see \cite[Theorem 2.10]{frediani2015shimura}).
For a cyclic group $G$ of size $m$, after identifying $\tg$ with $\Z/m\Z$, we have (see \cite[Lemma 2.7, Section 3.2]{moonen2010special})
\begin{align}
    g=1+\frac{(N-2)m-\sum_{j=1}^r\gcd(a(j),m)}{2},
\end{align}
and
\begin{align}\label{signature}
    f(\tau_i)=-1+\sum_{k=1}^r \langle\frac{-ia(k)}{m}\rangle \text{ for } 1\leq i\leq m-1, \text{ and } f(\tau_0)=0.
\end{align}

By Equation \eqref{signature}, for each orbit $\OO\in\mathfrak{O}_G$, the value $g(\tau)=f(\tau)+f(\tau^*)$ is independent of the choice of $\tau\in\OO$. In the following, we write  $g(\OO)=g(\tau)$, for any/all $\tau \in \OO$. Note that $g(\tau)=g(\tau^*)$\footnote{Throughout this paper, the assumption that $m$ is odd could be replaced  by the weaker assumption $f(\tau_{m/2})=0$ when $m$ is even}.

We make the following observations.
\begin{lemma}\label{lem: ft+fts<=r-2}
Let $(m,r,\au)$ be a cyclic monodromy datum.  Then,   for all $\tau \in \tg$, we have $
g(\tau)\leq r-2$.
\end{lemma}
\begin{proof}
By Equation \eqref{signature}, for every $\tau \in \tg$, we have the equation
\begin{equation}\label{eq: sum of conjugate signature}
g(\tau)=f(\tau)+f(\ts)
=-2
+\sum_{k=1}^r 
\Big(\langle\tfrac{-\tau a(k)}{m}\rangle 
+ \langle\tfrac{\tau a(k)}{m}\rangle \Big).
\end{equation}
Notice that 
$$
\Big\langle\frac{-\tau a(k)}{m}\Big\rangle 
+\Big\langle\frac{\tau a(k)}{m}\Big\rangle 
= \begin{cases}
    0 & \text{ if $m \mid \tau a(k)$} \\
    1 & \text{ if $m \nmid \tau a(k)$}.
\end{cases}
$$
Hence, $f(\tau)+f(\ts) \leq r-2$.
\end{proof}

\begin{lemma}\label{lem: inequality}
Let $(m,4,\au)$ be a cyclic monodromy datum with the number of branched points equal to $4$.  
Then,  there exists $\tn \in \mathcal{T}_{G}$ such that $f(\tn)=1$ and $f(\tns )=1$. 

Moreover, for all $\tau \in \OO_{\tn}$ and $1 \leq k \leq 4$, we have 
\begin{equation}
\begin{split}
\frac{p}{m}-1 
\leq \lfloor{p\langle \tfrac{\tau a(k)}{m}\rangle}\rfloor 
\leq p\frac{m-1}{m}.
\end{split}
\end{equation}
\end{lemma}
\begin{proof}
We prove the existence of of $\tn$ satisfying $f(\tn)=1$ and $f(\tns)=1$, and hence $g(\tn)=2$. Assume by contradiction that  $f(\tau) \neq 1$, for all $\tau \in \tg$.
By Lemma \ref{lem: ft+fts<=r-2}, for all $\tau \in \tg$, we have $f(\tau)+f(\ts) \leq 2$.
Hence, for each conjugate pair $(\tau, \ts)$, either $f(\tau)=0$ or $f(\ts)=0$. By construction,  the dimension of $\Sh(\D)$ is equal to $\sum_{\tau \in \tg}f(\tau)f(\ts)$. Hence, the assumption implies $\dim\Sh(\D)=0$  whereas $\dim(\M(G,4,\au))=1$, a contradiction. 

Let $\tn\in\tg$ satisfying $f(\tn)=1$ and $f(\tns)=1$. Then 
\[
f(\tn)+f(\tns)
=2
=-2
+\sum_{k=1}^4 
\Big(\langle\tfrac{-\tn a(k)}{m}\rangle 
+ \langle\tfrac{\tn a(k)}{m}\rangle \Big).
\]
We deduce that, for any $1 \leq k \leq 4$,  $m \nmid \tn a(k)$, and hence also $m \nmid p^i\tn a(k)$, for any $1 \leq i \leq l$. Hence, for any $\tau \in \OO_{\tn}$, we have
\[
\frac{1}{m} 
\leq \Big\langle \frac{\tau a(k)}{m}\Big\rangle 
\leq \frac{m-1}{m},
\]
which implies the statement.
\end{proof}

\subsection{The PEL-type moduli space $\Sh(\D)$}\label{Shimura}  
Given a monodromy datum $(G,r,\au)$, we define $S(G,r,\au)$ to be the largest closed, reduced and irreducible substack of $\mathcal{A}_{g,\C}$ containing the image of $\mgra_\C$ under the  Torelli morphism such that the multiplication of $\Z[G]$ on the Jacobian of the universal curve over $\mgra_\C$ extends to the universal abelian scheme over $S(G,r,\au)$.
By construction, $S(G,r,\au)$ is an irreducible connected component of the geometric fiber $\Sh(\D)_\C$ of a PEL-type  moduli space $\Sh(\D)$. The space $\Sh(\D)$ is a smooth proper Deligne-Mumford stack over $\OO_E[ \frac{1}{e}]$; furthermore, the Torelli morphism $\mgra_\C\to S(\gra)$ extends to a morphism
$\mgra \to \Sh(\D)$ over $\OO_E[\frac{1}{e}]$.
We denote by $\Sh(\D)_{\fpb}$ the geometric reduction of $\Sh(\D)$ modulo $p$, $p\nmid e$. 

In this paper, we consider two stratifications by locally closed strata of $\Sh(\D)_{\fpb}$, associated with two discrete invariants of abelian varieties: the  Newton polygon  and the Ekedahl-Oort type.

The \textbf{Newton polygon} is a discrete invariant that classifies the isogeny class of the $p$-divisible group of an abelian variety over $\fpb$. By \cite{viehmann2013ekedahl}, the Newton polygons corresponding to non-empty strata in $\Sh(\D)_{\fpb}$ are in one-to-one correspondence with the elements in the partially ordered Kottwitz's set $B_p(\D)$ (see \cite{kottwitz1985isocrystals}). By \cite[Theorem 1.6.2 (Density Theorem)]{wedhorn1999ordinariness}, the unique maximal element in $B_p(\D)$ corresponds to the unique non-empty and dense Newton polygon stratum in $\sh(\D)_{\fpb}$, called the  $\mu$-ordinary locus. 

The \textbf{Ekedahl-Oort type} is a discrete invariant that classifies the isomorphism class of the $p$-kernel of a polarized abelian variety over $\fpb$, together with its multiplication by $\F_p[G]$. By \cite[Theorem 2 (1),(2)]{viehmann2013ekedahl}, the Ekedahl-Oort types corresponding to non-empty strata in $\shgf_{\fpb}$ are in one-to-one correspondence with certain elements in the Weyl group ${\rm Weyl}(\D)$ of the reductive group associated with the datum $\D$; furthermore,  the unique element of maximal length  corresponds to the unique non-empty open and dense Ekedahl-Oort stratum in  $\shgf_{\fpb}$, called the $p$-ordinary locus.

By \cite[Theorem 1.3.7]{moonen2004serre}, the $p$-ordinary Ekedahl-Oort stratum and the
$\mu$-ordinary Newton polygon stratum of $\sh(\D)_{\fpb}$ agree. 
By \cite{Goldring2013TheH}, this stratum is the non-vanishing locus of an automorphic form on $\Sh(\D)_{\fpb}$, called the \textbf{$\mu$-ordinary Hasse invariant}.

Under the Torelli morphism, the stratifications by Newton polygons and \EO-types on $\sh(\D)_{\fpb}$ induces stratifications on the Hurwitz space $\M(G,r,\au)$. As their definition depends on the PEL-datum $\D$ and the prime $p$, we refer to the pullback on $\mgra_{\fpb}$ of the $\mu$-ordinary/$p$-ordinary stratum  as the $\D$-ordinary locus.

\begin{definition}[$\D$-ordinary Hasse--Witt invariant]\label{def: D ordinary Hasse invariant}
Given a monodromy datum $(G,r,\au)$,  a $\D$-ordinary Hasse--Witt invariant at a prime $p$, $p\nmid |G|$, is a polynomial $h \in \fpb[x_1, \dots, x_r]$ whose non-vanishing locus is equal to the $\D$-ordinary locus of $\mathcal{M}(G,r,\au)_{\fpb}$.   
\end{definition}
In \cite{lin2024abelian}, under the assumption that $r \leq 5$ and $G$ is cyclic, we compute a $\D$-ordinary Hasse--Witt invariant on $\mgra_{\fpb}$ via the notion of Hasse--Witt triple.

\subsection{The Hasse--Witt triple of cyclic covers of $\Po$}
In \cite{moonen2022computing}, Moonen establishes a equivalence of category between the polarized mod-$p$ Dieudonn\'e modules and Hasse--Witt triples; by definition, a Hasse-Witt triple $(Q,\phi,\psi)$ is as follows:
\begin{itemize}
    \item $Q$ is a finite dimensional vector space over $\fpb$;
    \item $\phi: Q \to Q$ is a $\sigma$-linear map;
    \item $\psi: \ker(\phi) \to \im(\phi)^{\perp}$ is a $\sigma$-linear isomorphism, where $\im(\phi)^{\perp}$ is the subspace of $Q^{\smvee}$ such that $\im(\phi)^{\perp}=\{\lambda \in Q^{\smvee}\mid \forall q \in Q: \lambda(\phi(q))=0 \}$.
\end{itemize}
In \cite{moonen2022computing}, Moonen gives an explicit algorithm for computing the Hasse--Witt triple of a complete intersection curve defined over $\fpb$, where $Q=\hy$, $Q^{\smvee}=\hyd$, and $\phi:\hy \to \hy$ is induced by the action of Frobenius (the Hasse--Witt matrix).
In \cite{lin2024abelian}, we adapt Moonen's algorithm to the context of cyclic covers of $\Po$. We briefly recall the computation. 

Let $\ux=(x_1,\dots,x_r)\in\mlra(\fpb)$ and denote $C=C_{\ux}$. 
The action of $G$ on $Q$  yields a decomposition $Q=\bigoplus_{\tau\in \tg} Q_{\tau}$, where $Q_{\tau}$ is the $\tau$-isotypic component of $Q$, with $\dim(Q_{\tau})=f(\tau^*)$ (here, $\fb$ denotes the signature of $C$,  computed in \eqref{signature}). Correspondingly, there is a decomposition on the dual space, $Q^{\vee}=\bigoplus_{\tau\in \tg} Q_{\tau}^{\vee}$ (here $Q_{\tau}^{\vee}=(Q_{\tau})^{\vee}\neq (Q^{\vee})_{\tau}$).
Let $\phi_{\tau}$ and $\psit$ be the restrictions of $\phi$ and $\psi$ respectively to $Q_{\tau}$. Then, $\phit:Q_{\tau}\to Q_{p\tau}$ and $\psit:Q_{\tau}\to Q_{p\tau}^{\vee}$. 

In \cite[Section 5]{lin2024abelian},  for each $1\leq i\leq m-1$, we describe explicit bases of $Q_{\tau_i}$ and $Q_{\tau_i}^{\vee}$ in terms of Cech cohomology.  With respect to those bases,  the map $\phi_{\tau_i}:Q_{\tau_i}\to Q_{p\tau_i}$ is computed as a $f(p\tau_i^*)\times f(\tau_i^*)$ matrix (called the Hasse--Witt matrix), whose $(j',j)$-entry is
\begin{align}\label{coefficient of phii}
\phi_{\tau_i}(j',j)
=(-1)^N 
\sum_{n_1+ \dots + n_r=N}
\tbinom{\pby}{n_1} \dots 
\tbinom{\pbr}{n_r} 
x_1^{n_1} \dots x_r^{n_r},    
\end{align} 
where $N=\sum_{k=1}^r \pbk -(jp-j')=p(f({\tau_i}^*)+1)-(f(p{\tau_i}^*)+1)-jp+j'$.

Under the additional assumptions that $\gcd(i,m)=1$ and either $f(p\tau_i^*)=0$ or
$f(\tau_i^*)=g(\tau_i)$,  we compute the map $\psi_{\tau_i}:Q_{\tau_i}\to Q_{p\tau_i}^{\vee}$ as a $f(p\tau_i)\times f(\tau_i^*)$ matrix whose $(j',j)$-entry is
\begin{align}\label{coefficient of psii}
\psi_{\tau_i}(j',j)
=-\sum_{k=1}^r 
{\lfloor{p\langle\tfrac{ia(k)}{m}\rangle}\rfloor}
r_{i,j,k}q_{r-j',k},
\end{align}
where
\begin{equation*}
    q_{r-j',k}  \text{ is the coefficient of } x^{j'-1} \text { in } \frac{(x-x_1) \dots (x-x_r)}{(x-x_k)},
\end{equation*}

$$
r_{i,j,k}
=(-1)^N
\sum_{n_1+ \dots +n_r=N}
\tbinom{\pby}{n_1} 
\dots 
\tbinom{{\lfloor{p\langle\frac{ia(k)}{m}\rangle}\rfloor} -1}{n_k} 
\dots \tbinom{\pbr}{n_r} 
x_1^{n_1}\dots x_r^{n_r},
$$
and $N=\sum_{k=1}^r {\lfloor{p\langle\tfrac{ia(k)}{m}\rangle}\rfloor} -jp=p(f(\tau_i^*)+1)-(f(p\tau_i^*)+1)-pj$. 

In the following, we denote by $\check{\phi}_\tau: Q_{\tau}^{\vee} \to Q_{p\tau}^{\vee}$  (resp. $\check{\psi}_\tau: Q_{\tau}^{\vee} \to Q_{p\tau}$) the maps induced via duality  by $\phi_\tau$ (resp.  $\psi_\tau$). With respect to our choice of bases of $Q_\tau$, $Q_\tau^{\smvee}$, the matrices representing $\phi_\tau$ and $\check{\phi}_\tau$ (resp. $\psi_\tau$ and $\check{\psi}_\tau$) are the same.

\subsection{$\D$-ordinary $\OO$-partial Hasse--Witt invariant}\label{Sec: mu ord hasse inv}

In this subsection, we recall the computation of the $\D$-ordinary Hasse--Witt invariant from \cite{lin2024abelian}, and adapt it to our current context.
\begin{notation}\label{notation: notation 0}
Let $(m,r,\au)$ be a cyclic monodromy datum, with $r\geq 4$, and $p$ a prime, $p\nmid m$. 
Assume (A1): there exists $\tau \in \tg$ where $f(\tau)=1$, and $g(\tau)\geq 2$.  

We fix  $\tn\in \tg$ satisfying $f(\tns)=1$, and $g(\tau_0^*)\geq 2$. Denote
 $\OO =\OO_{\tn}$  the Frobenius orbit of $\tn$, and $l=l(\OO)$ the size of $\OO$.  We identify  the set $\{0, 1, \dots, l-1\}$ with $\OO$  via $i\mapsto  p^i\tn $.
\end{notation}

Following the notations in \cite[Section 7]{lin2024abelian}, for each $\tau \in \OO$, we write
\begin{equation}\label{eq: Eq_A_i}
    A_{\tau}=
    \begin{cases}
    \phi_{\tau}: Q_\tau\to Q_{p\tau} &\text{if } f(\ts)\geq 1, f(p\ts)\geq 1,\\
    \check{\phi}_{\tau}: Q_{\ts}^{\smvee} \to Q_{p\ts}^{\smvee}  &\text{if } f(\ts)=0, f(p\ts)=0,\\
    \psi_{\tau}: Q_{\tau}\to Q_{p\ts}^{\smvee} &\text{if } f(\ts)\geq 1, f(p\ts)=0,\\
    \check{\psi}_{\tau}: Q_{\ts}^{\smvee} \to Q_{p\ts} &\text{if } f(\ts)=0, f(p\ts)\geq 1.
    \end{cases}
\end{equation}
With respect to our choice of bases,  $A_\tau$ is a matrix of size $d(p\tau) \times d(\tau)$, where
\begin{equation}\label{def: ai}
    d(\tau)=
    \begin{cases}
      f(\ts) &\text{ if } f(\ts)\geq 1,\\
      f(\tau)  & \text{ otherwise.}
    \end{cases}
\end{equation}
The entries of the matrix $A_\tau$ are homogeneous polynomials in the variables $x_1, \dots,x_r$. 
In the following, we write $A_i=A_{p^i\tn}$ and $d(i)=d(p^i\tn)$, for $0 \leq i \leq l-1$.

\begin{definition}\label{def:h_0}

(\cite[Equation (20)]{lin2024abelian})
With notations as in Notation \ref{notation: notation 0},  we define the $\D$-ordinary $\OO$-partial Hasse--Witt invariant $h_0$ as 
 \begin{equation}\label{h_0}
h_0 =A_{l-1}\circ \cdots A_1\circ A_0  \in\fp[x_1,\dots,x_r]. \end{equation}
\end{definition}

By \cite[Theorem 7.2]{lin2024abelian},  
 the vanishing locus of $h_0$ {is contained in } the non-$\D$-ordinary locus of the compactified Hurwitz space $\overline{\M(G,r,\au)}$. Furthermore, when $r\leq 5$, by \cite[Theorem 7.2]{lin2024abelian}, the non-$\D$-ordinary locus is the union of the vanishing loci of the $\D$-ordinary $\OO$-partial Hasse--Witt invariants from Defintion \ref{def:h_0}, as $\OO$ varies among the Frobenius orbits in $\tg$ (see Remark \ref{partialhasse}).
 In particular, we deduce the following statement.
 
 \begin{lemma}\label{lemma:reduction}
 (\cite[Theorem 7.2]{lin2024abelian})
If $\alpha=(\alpha_1,\dots, \alpha_r)\in \fpb^r$ satisfies $h_0(\alpha)=0$ and $\alpha_1, \dots, \alpha_r\in\fpb$ are all distinct, then 
$\alpha$ defines a non-$\D$-ordinary point in $\M(m,r,\au)(\fpb)$.     
 \end{lemma}
 
In the following, we establish the existence of smooth non-$\D$-ordinary covers over $\fpb$, by proving that the vanishing locus of $h_0$ is not contained in the boundary of $\overline{\M(G,r,\au)}_{\fpb}$. 


We observe that by \cite[Equations (19) and (20)]{lin2024abelian}, the $\D$-ordinary $\OO$-partial Hasse--Witt invariant $h_0$  in Equation \eqref{h_0} satisfies 
\begin{equation}\label{RJi}
    h_0=\sum_{J\in\mathfrak{J}} R_{J} \text{  with  } R_J=\prod_{i=0}^{l-1}R_{J,i}^{p^{l-i-1}} ,
\end{equation}
where \begin{itemize}
    \item 
$\mathfrak{J}$ is the set of all functions $J:\{0,1, \dots, l\} \to \mathbb{N}$
satisfying $J(0)=J(l)=1$, and  $1 \leq J(i) \leq d(i)$  for each $0<i<l$;
\item for $J \in \mathfrak{J}$ and $0<i<l$,  
 $R_{J, i}=A_i(J(i+1),J(i))\in\fp[x_1,\dots,x_r]$ is the $(J(i+1),J(i))$-entry of the matrix $A_i$ from Equation \eqref{eq: Eq_A_i}. 
\end{itemize}
\begin{lemma}\label{degree}
With notations as in Notation \ref{notation: notation 0}, let $h_0$ be the $\D$-ordinary $\OO$-partial Hasse--Witt invariant defined by Equation \eqref{h_0}. 
    
    Then, 
    $\deg(h_0)>p^{l-1}(p-(r-2))$.
\end{lemma}
\begin{proof} Define $J_0\in\mathfrak{J}$ by $J_0(i)=1$, for all $0 \leq i \leq l-1$.  By \cite[Proposition 6.3]{lin2024abelian}, $R_{J_0}\neq 0$. Hence,
    from Equation \eqref{RJi}, we deduce 
\begin{align*}
\deg(h_0)&=\max_{J}(\deg R_J) \geq \deg R_{J_0} > p^{l-1}\deg(R_{J_0,0}).
\end{align*}
Furthermore, the assumption $f(\tau_0^*)=1$ implies
$\deg(R_{J_0,0})\in\{ {s}-p+1, s-p+r-1\}$ by Equation \eqref{eq: Eq_A_i}.
In particular, $\deg(R_{J_0,0})\geq s-p+1=p-f(p\tau^*)\geq p-(r-2)$.
\end{proof}

\subsection{Reduction to the cyclic case}\label{sec: reduction}

Following \cite[Section 3]{lin2024abelian}, we reduce the proof of Theorems \ref{main theorem} and \ref{prop: main theorem for r eq 4} to the special case of cyclic monodromy data, Propositions \ref{main prop} and \ref{prop:main r=4} respectively. 

\begin{proof}[Proof of Theorems \ref{main theorem} and \ref{prop: main theorem for r eq 4}]
Let $\tau_0 \in \tg$ satisfying $f(\tau_0)=1$ and $g(\tau_0)\geq 2$. 
Set $Q=G/\ker(\tau_0)$. Then $Q$ is cyclic, of size $|Q|\leq |G|$, and $\tau_0$ descends to a character $\bar{\tau}_0\in \mathcal T_Q$.   
By Lemma 1, when $r=4$, 
if $g(\tau_0)\geq 2$ then $g(\tau_0)=2$; moreover, for any prime $p\nmid |G|$,  $\ker (p\tau_0)=\ker(\tau_0)$ and hence $p\tau_0$ also descends to $Q$.

Let $C$ be a smooth project $G$-cover of $\Po$, with monodromy datum $(G,r,\au)$, and signature type $\fb$, and denote by $\D$ the associated PEL-datum. By construction, the quotient curve $C/\ker(\tau_0)$ is a $Q$-cover of $\Po$. We denote by $(Q,r',\ub)$ its monodromy datum, and by $\fb_Q$ its signature type. Then, by construction, $f_Q(\bar{\tau}_0)=f(\tau_0)=1$, and $g_Q(\bar{\tau}_0)=g(\tau_0) \geq 2$. Moreover, since $\bar{\tau}_0\in\mathcal{T}_Q$ has trivial kernel,  $g_Q(\bar{\tau}_0)=r'-2$, hence  $r' \geq 4$.

Let $\D_Q$ denote the PEL-datum associated to the monodromy datum $(Q, r',\ub)$. 
By Lemma \cite[Lemma 3.2]{lin2024abelian}, if
curve $C/\ker(\tau_0)$ is non-$\D_Q$-ordinary, then $C$ is non-$\D$-ordinary. 
Since the projection $\mgra\to \M(Q,r',b)$, mapping the isomorphism class of $C$ to that of $C/\ker(\tau_0)$, is surjective on geometric points, of characteristic $p\nmid |G|$, the statements follow from Propositions \ref{main prop}\ and \ref{prop:main r=4}.
\end{proof}

\section{The extended Hasse--Witt matrix of cyclic covers branched at four points}\label{sec: hasse invariant one dim}

In this section, we focus on cyclic covers of $\Po$ branched at four points ($r=4)$, and explicitly compute the entries of the matrices $A_\tau$ from Equation \eqref{eq: Eq_A_i} appearing in the definition of the $\D$-ordinary $\OO$-partial Hasse--Witt invariant $h_0$ from Definition \ref{def:h_0}. 


We recall notations and assumptions from Section \ref{Sec: mu ord hasse inv}

\begin{notation}\label{notation: notation one}
We fix a cyclic monodromy datum $(m,4,\au)$ with 4 branched points,  a prime $p>3m$,  and a character $\tn\in\tg$ satisfying $f(\tn)=1$ and $f(\tns)=1$. 

After reordering the entries of $\au$, we assume
\begin{equation}\label{order}
\floor*{p\left\langle\tfrac{\tn a(2)}{m}\right\rangle}
+
\floor*{p\left\langle\tfrac{\tn a(3)}{m}\right\rangle}
\leq
\floor*{p\left\langle\tfrac{\tn a(1)}{m}\right\rangle}
+
\floor*{p\left\langle\tfrac{\tn a(4)}{m}\right\rangle}.
\end{equation}

For any $\tau\in \OO=\OO_{\tn}$ and integer $N$, we define
\begin{align}
\st&=\sum_{k=1}^4 {\floor*{p\langle\tfrac{\tau a(k)}{m}\rangle}}=p(f({\tau}^*)+1)-(f(p{\tau}^*)+1), \label{st} \\
c_{\tau}(N)
&=\min\left\{1 \leq q \leq 4: \tsum_{k=1}^q {\floor*{p\langle \tfrac{\tau a(k)}{m}\rangle}}
>N\right\}.\label{ct}
\end{align}

\end{notation}


\begin{lemma}\label{c23}\label{Lem: st - two floors}
Suppose $f(\tau^*)\neq 0$.
Then, we have \begin{enumerate}
    \item $ 3p-3\leq s_\tau\leq 3p-1$  if $f(\ts)=2$, and $2p-3\leq s_\tau\leq 2p-1$ if $f(\ts)=1$.

\item  
$
2\leq c_{\tau}(s_{\tau}-pj+j')
\leq 3,
$ for any $1 \leq j \leq f(\ts)$ and $0 \leq j' \leq 2$.
\item 
$p<\st 
-\lfloor{p\tfrac{\tau a(k_1)}{m}}\rfloor
-\lfloor{p\tfrac{\tau a(k_2)}{m}}\rfloor
<2p-2,
$  for any pair $1\leq k_1<k_2\leq 4$, if $\fts=2$.
\end{enumerate}
\end{lemma}

\begin{proof}
The bounds on $s_\tau$ follow from the definition, while
by definition, we have $1\leq c_\tau(s_{\tau}-pj+j')\leq 4$.
By \cite[Lemma 6.2]{lin2024abelian}, $c_\tau(s_{\tau}-pj+j')\neq 4$. 
By assumption, 
\[s_\tau-pj+j'= p(\fts+1-j)-(f(p\ts)+1-j')\geq p-3.\]
Since $\pbtoinl\leq p\frac{m-1}{m}<p-3$,  we deduce
 $\pbtoinl < s_\tau-pj+j'$, and hence $c_\tau(s_{\tau}-pj+j')\neq 1$. 

Finally, let $\{k_1,k_2\}\subset \{1,2,3,4\}$ and set $\{k_3,k_4\}=\{1,2,3,4\}\setminus \{k_1,k_2\}$. 
Since $\fts=2$, we have
$\st
\geq 3p-3$; since $p>3m$, we deduce
\[
\st
-\lfloor{p\tfrac{\tau a(k_1)}{m}}\rfloor
-\lfloor{p\tfrac{\tau a(k_2)}{m}}\rfloor 
\geq 3p-3
-p\left(\tfrac{m-1}{m}\right)
-p\left(\tfrac{m-1}{m}\right)
 \geq p-3+ \frac{2p}{m}
>p,
\]
 and
\[
\st 
-\lfloor{p\tfrac{\tau a(k_1)}{m}}\rfloor
-\lfloor{p\tfrac{\tau a(k_2)}{m}}\rfloor
=\lfloor{p\tfrac{\tau a(k_3)}{m}}\rfloor
+\lfloor{p\tfrac{\tau a(k_4)}{m}}\rfloor
\leq 2 p \frac{m-1}{m}
<2p-2. \qedhere
\]
\end{proof}
By Lemma \ref{c23}, since $f(\tns)=1$, we have
$s_{\tn}
=
\sum_{k=1}^4
\floor*{p\left\langle\tfrac{\tn a(k)}{m}\right\rangle}
\leq 2p-1.$
Together with Equation \eqref{order}, this implies
$\floor*{p\left\langle\tfrac{\tn a(2)}{m}\right\rangle}
+
\floor*{p\left\langle\tfrac{\tn a(3)}{m}\right\rangle}
\leq p-1$.

When $r=4$,  by a linear fractional transformation on $\Po$, without loss of generality, we may assume the branched points to have coordinates $(x_1,x_2,x_3,x_4)=(\infty, t, 1,0)$. 

\begin{notation}\label{notation:t}
    For any  polynomial $F=F(x_1, x_2, x_3, x_4) \in \fp[x_1, x_2, x_3, x_4]$, we denote by $F(t)\in \fp[t]$ its  specialization to  $(x_1,x_2,x_3,x_4)=(\infty, t, 1,0)$. 
\end{notation}

In Proposition \ref{specialization of phi and psi},  we compute the entries of $\phi_\tau$ and $\psi_\tau$ at $(x_1,x_2,x_3,x_4)=(\infty, t, 1,0)$ in terms of the following class of polynomials (which we study in  Appendix~\ref{sec: separability}). 
\begin{definition}\label{def: fabc}
For any non-negative integers $a,b,c$, we define $f(a,b,c)\in \fpb[t]$ as  
\begin{equation}\label{fabc}
  f(a,b,c)=\sum_{i_2+i_3=c}\tbinom{a}{i_2}\tbinom{b}{i_3}t^{i_2}. 
\end{equation}
\end{definition}

{We summarize the properties of polynomials $f(a,b,c)$ proven in Appendix \ref{sec: separability}.}
\begin{lemma}\label{Lem: deg, vt and vt-1 of fabc}
Let $a,b,c\geq 0$ be integers. Assume $a,b<p$ and $c\leq a+b$.

Then, the polynomial $f(a,b,c)\in \fp[t]$  satisfies the following properties
\begin{enumerate}
\item  it has only simple roots except possibly $0$ and $1$.
  
\item it has value at $1$ equal to $f(a,b,c)(1)=\binom{a+b}{c}$
\item it vanishes  at 1 with multiplicity $$v_{t-1}(f(a,b,c)) = 
a+b-(p-1)    \text { if }
  0< a+b-(p-1)\leq c\leq p-1  , \text{ and 0 otherwise.}$$
  \item it has degree equal to $\deg_t(f(a,b,c))=\min\{a,c\}$;
\item it vanishes at $0$ with multiplicity  $ v_t(f(a,b,c)) =\max\{0,c-b\}$;
\end{enumerate}
\end{lemma}

\begin{proposition}\label{specialization of phi and psi}
The entries  of $\phi_{\tau}$ and $\psi_{\tau}$ from Equations \eqref{coefficient of phii} and \eqref{coefficient of psii} satisfy
\[\phi_\tau(j',j)(t)
=(-1)^{\st-pj+j'}
f\left(\pbtt, \pbtd, \st-pj+j'-\pbto\right);\]
$\psit(1,j)(t)
=(-1)^{\st-pj+2}\pbtf t
f\left(\pbtt, \pbtd, \st-pj-\pbto  \right)$;
\[\begin{split}
&\psit(2,j)(t)=\\
&(-1)^{\st-pj+2}
\left(\st-\pbto+1\right)
f\left(\pbtt,\pbtd, \st-pj-\pbto+1\right).
\end{split}\]
\end{proposition}

\begin{proof}
{Note that, given a polynomial $F \in \fp[x_1,x_2,x_3,x_4]$, the specialization $F(t)$ of $F$ at $(x_1,x_2,x_3,x_4)=(\infty,t,1,0)$ is equal to the sum of all monomials $M$ in $F$, for which the difference $v_{x_1}(M)-v_{x_4}(M)$ is maximal, specialized at  $(x_2,x_3)=(t,1)$}

Let $\phit(j',j)$ as defined in Equation \eqref{coefficient of phii}; by Lemma \ref{c23}, the monomials $M$ in $\phit(j',j)$ for which the difference $v_{x_1}(M)-v_{x_4}(M)$ is maximal satisfy $v_{x_1}(M)=\pbto$ and $v_{x_4}(M)=0$. Summing over all such monomials gives the stated formula for $\phit(j',j)$.

Let $\psit(j',j)$ as defined in Equation \eqref{coefficient of psii};  any monomial $M$ in $\psitjjp$ can be factored as $M=M_1M_2$, where  $M_1$ are the monomials in $r_{\tau,j,k}$ and $M_2$ are the monomials in $q_{4-j',k}$, for some $1\leq k\leq 4$. By Lemma \ref{c23}, the difference $v_{x_1}(M)-v_{x_4}(M)$ is maximal if the monomial M satisfies $v_{x_1}(M)=\pbto+1$ and $v_{x_4}(M)=0$; that is, when $M=M_1M_2$ with $M_1,M_2$ satisfying  $v_{x_1}(M_1)=\pbto$, $v_{x_4}(M_1)=0$, and  $v_{x_1}(M_2)=1$, $v_{x_4}(M_2)=0$.

When $j'=1$, we have $q_{3,k}=-\frac{x_1x_2x_3x_4}{x_k}$, $1\leq k\leq 4$; hence, $M_2=q_{3,k}$ and it satisfies $v_{x_1}(M_2)=1$ and $v_{x_4}(M_2)=0$ if and only if $k=4$. We deduce the stated formula for
\(\psi_\tau(1,j)(t)\).


   When $j'=2$, then there exists a monomial $M_2$ in $q_{2,k}$ with $v_{x_1}(M_2)=1$ and $v_{x_4}(M_2)=0$ if and only if $k\neq 1$. Furthermore, $M_2$ is equal to $x_1x_3$ if $k=2$, to $x_1x_2$ if $k=3$, and to either $x_1x_2$ or $x_1x_3$ if $k=4$. We deduce
    \begin{align*}
    \psit(2,j)(t)
    &=(-1)^{\st-pj+3}\pbtt
    \sum_{i_2+i_3=\st-pj-\pbtoinl}\tbinom{\pbttinl-1}{i_2}
    \tbinom{\pbtdinl}{i_3}t^{i_2} \\
    &+(-1)^{\st-pj+3}
    \pbtd t
    \sum_{i_2+i_3=\st-pj-\pbtoinl}
    \tbinom{\pbttinl}{i_2}
    \tbinom{\pbtdinl-1}{i_3}
    t^{i_2}\\
    &+(-1)^{\st-pj+3}
    \pbtf(t+1)
    \sum_{i_2+i_3=\st-pj-\pbtoinl}
    \tbinom{\pbttinl}{i_2}\tbinom{\pbtdinl}{i_3}t^{i_2}.
    \end{align*}
    Set $a=\st-pj-\pbto$.
    Ignoring the sign $(-1)^{\st-pj+3}$, the coefficient of $t^i$ in $\psit(2,j)(t)$ is 
    \begin{align*}
    &\pbtt\tbinom{\pbttinl-1}{i}
    \tbinom{\pbtdinl}{a-i} 
    +\pbtd
    \tbinom{\pbttinl}{i-1}
    \tbinom{\pbtdinl-1}{a-i+1}\\
    &+\pbtf\tbinom{\pbttinl}{i-1}
    \tbinom{\pbtdinl}{a-i+1} 
    +\pbtf\tbinom{\pbttinl}{i}
    \tbinom{\pbtdinl}{a-i}\\
    &=\tbinom{\pbttinl}{i}
    \tbinom{\pbtdinl}{a-i+1}
    \cdot
    \Big(\pbtt
    \frac{\pbtt-i}{\pbtt} \frac{a-i+1}{\pbtd-(a-i)} \\
    &+\pbtd\frac{i}{\pbtt-i+1} \frac{\pbtd-a+i-1}{\pbtd} \\
    &+\pbtf \frac{i}{\pbtt-i+1} 
    +\pbtf\frac{a-i+1}{\pbtd-a+i} \Big),
    \end{align*}
which equals $-(a+1)\tbinom{\pbttinl}{i}\tbinom{\pbtdinl}{a-i+1}$, where $a\equiv\pbtt+\pbtd+\pbtf\pmod{p}$.  We deduce the stated formula for
$\psi_\tau(2,j)(t)$, after multiplying back by the sign $(-1)^{s_\tau-pj+3}$.
\end{proof}

\begin{lemma}\label{verify}

   The polynomials $f(a,b,c)$ occurring in Proposition \ref{specialization of phi and psi} satisfy the assumptions of Lemma \ref{Lem: deg, vt and vt-1 of fabc}, namely  $a,b<p$ and $c< a+b$.
\end{lemma}
\begin{proof}
First,  $a,b<p$ since $a=\pbtt$, and $b=\pbtd$. 

\begin{enumerate}
    \item For $\phi_\tau(j',j)(t) $: 
    \begin{align*}
    &a+b-c=
\pbtt+\pbtd-(\st-pj+j'-\pbto) \\
&= pj-j'-\pbtf > p-2-p\frac{m-1}{m} >0.
\end{align*}
\item For $\psi_\tau(1,j)(t)$: $a+b-c=pj-\pbtf>0$
\item For $\psi_\tau(2,j)(t)$: $a+b-c=pj-1-\pbtf>0$. \qedhere
\end{enumerate}
\end{proof}

By Proposition \ref{specialization of phi and psi} and Lemma \ref{verify}, we  specialize Lemma \ref{Lem: deg, vt and vt-1 of fabc} to $\phi_\tau(j',j)(t)$ and $\psi_\tau(j',j)(t)$.

 \begin{corollary}
[Separability outside of $0$ and $1$]\label{prop: separability of phit and psit}
For $1\leq j'\leq f(p\ts)$,  $1\leq j\leq f(\tau^*)$, the polynomial
$\phit(j',j)(t)$ and $\psit(j',j)(t)$ have only simple roots, except possibly $0$ and $1$.
\end{corollary}

\begin{corollary}
[Value at $1$]\label{Lem: value at one}
Let $1\leq j'\leq 2$,  $1\leq j\leq f(\tau^*)$. Then, 
\begin{equation}\label{Eq: phit at one}
\phi_{\tau}(j',j)(1)
=(-1)^{\st-pj+j'} \binom{\pbttinl+\pbtdinl}{\st-pj+j'-\pbtoinl}
\end{equation}
\begin{equation}\label{Eq: psit at one}
    \psi_{\tau}(j',j)(1)
=(-1)^{\st-pj+j'+1}\pbtf \binom{\pbttinl+\pbtdinl}{\st-pj-\pbtoinl}
\end{equation}
\end{corollary}
\begin{proof}

For $\phi_{\tau}(j',j)$ and $\psi_{\tau}(1,j)$ the statement follow immediately; for $\psi_{\tau}(2,j)$, we have
\begin{align*}
&\psi_{\tau}(2,j)(1)
=(-1)^{\st-pj+2}(\st-\pbto+1) 
\binom{\pbttinl+\pbtdinl}{\st-pj-\pbtoinl+1}\\
&=(-1)^{\st-pj+2}\Big(\pbtt+\pbtd-(\st-pj-\pbto)\Big) 
\binom{\pbttinl+\pbtdinl}{\st-pj-\pbtoinl}\\
&=(-1)^{\st-pj+2} \Big(-\pbtf \Big) 
\binom{\pbttinl+\pbtdinl}{\st-pj-\pbtoinl}.\qedhere
\end{align*}
\end{proof}

\begin{corollary}[Vanishing at 1]\label{Cor: mult of one}
Let $1\leq j'\leq f(p\ts)$,  $1\leq j\leq f(\tau^*)$. Assume either
$\pbttinl+\pbtdinl\leq p-1$
or $f(\ts)=2$ and $j=1$.
Then, the polynomials $\phi_{\tau}(j',j)$ and $\psi_{\tau}(j',j)$ do not vanish at $1$.
\end{corollary}
\begin{proof}
It suffices to note that if  $f(\ts)=2$ and $j=1$, then $s_{\tau}\geq 3p-3$; hence, 
$
\st-p-\pbto+1\geq 2p-3-\pbto
\geq 2p-3- p \tfrac{m-1}{m}
\geq p$ since $p>3m$.
\end{proof}

\begin{corollary}\label{Cor: deg}
Suppose $f(\ts)=2$.
Then, the polynomials $\phit(j',j)(t)$ and $\psit(j',j)(t)$ have degree

\begin{enumerate}
    \item $\deg_{t}(\phit(j',j)(t))=\begin{cases}
        \pbtt & \text{ if $j=1$, }\\
        \st-pj+j'-\pbto & \text{ if $j=2$}
    \end{cases}$
\item
$\deg_{t}(\psit(1,j)(t))=\begin{cases}
1+\pbtt & \text{ if $j=1$,} \\
1+\st-pj-\pbto & \text{ if $j=2$}
\end{cases}$
\item $\deg_{t}(\psit(2,j)(t))=\begin{cases}
    \pbtt & \text{ if $j=1$,} \\
    \st-pj+1-\pbto & \text{ if $j=2$}
\end{cases}$
\end{enumerate}
\end{corollary}
\begin{proof}
The result follows from Proposition \ref{specialization of phi and psi} and Lemma \ref{Lem: deg, vt and vt-1 of fabc}(4), using Lemma \ref{c23}(3) to determine the minimum in each case.
\end{proof}

\begin{corollary}
[Vanishing at $0$]\label{Cor:vt}
Suppose $f(\ts)=2$. Then, the polynomials $\phi_\tau(j',j)$ and $\psi_\tau(j',j)$ vanish at $0$ with multiplicities
\begin{enumerate}
\item $v_{t}(\phi_{\tau}(j',j)(t)) =\begin{cases}
    \st-pj+j'-\pbto-\pbtd &\text{ if $j=1$ }, \\
    0 & \text{ if $j=2$ }.
    \end{cases}$
\item $v_{t}(\psi_{\tau}(1,j)(t)) = \begin{cases}
        \st-pj+1-\pbto-\pbtd &\text{ if $j=1$, } \\
        1 &\text{ if $j=2$}.
    \end{cases}$
\item $v_{t}(\psi_{\tau}(2,j)(t)) = \begin{cases}
        \st-pj+1-\pbto-\pbtd &\text{ if $j=1$, } \\
        0 &\text{ if $j=2$}.
    \end{cases}$
\end{enumerate}
\end{corollary}
\begin{proof}
The result follows from Proposition \ref{specialization of phi and psi} and Lemma \ref{Lem: deg, vt and vt-1 of fabc}(5), using Lemma \ref{c23}(3) to determine the maximum in each case.
\end{proof}

\begin{corollary}[Vanishing at $\infty$]\label{Cor: bound on degree}
Suppose $f(\ts)=2$.  Then,
\begin{enumerate}
    \item $\deg(\phit(j',j)(t))-v_t(\phit(j',j)(t)) \geq \frac{p}{m}-3$
    \item $\deg(\psit(j',j)(t))-v_t(\psit(j',j)(t)) \geq \frac{p}{m}-3$
\end{enumerate}
\end{corollary}
\begin{proof}
By Corollaries  \ref{Cor: deg} and \ref{Cor:vt}, we compute $\deg -v_t$,  the order of vanishing at $\infty$.

\begin{enumerate}
\item   $\deg(\phit(j',j)(t))-v_t(\phit(j',j)(t))=$ 
\begin{align*}
\begin{cases}
\pbtt-(\st-pj+j'-\pbto-\pbtd) & \text{ if $j=1$} \\
\st-pj+j'-\pbto & \text{ if $j=2$}.
\end{cases}
\end{align*}

\item 
  $\deg(\psit(1,j)(t))-v_t(\phit(2,j)(t))=$ 
\begin{align*}
\begin{cases}
\pbtt-(\st-pj-\pbto-\pbtd) & \text{ if $j=1$} \\
\st-pj-\pbto & \text{ if $j=2$}.
\end{cases}
\end{align*}

\item 
 $\deg(\psit(2,j)(t))-v_t(\psit(2,j)(t))=$ 
\begin{align*}
\begin{cases}
\pbtt-(\st-pj+1-\pbto-\pbtd) & \text{ if $j=1$} \\
\st-pj+1-\pbto & \text{ if $j=2$}.
\end{cases}
\end{align*}
\end{enumerate}

When $j=1$, all three values of $\deg -v_t$ are at least $p-\pbtfinl-2$, and

\[
p-\pbtf-2
\geq p-p\frac{m-1}{m} -2= \frac{p}{m}-2.
\]

When $j=2$,  all three values of $\deg -v_t$ are at least $\st-pj-\pbtoinl$, and
\[
\st-pj-\pbto  
\geq (3p-3) +2p -p\frac{m-1}{m}
= \frac{p}{m}-3.\qedhere
\]
\end{proof}

\section{Non-$\D$-ordinary covers branched at four points}\label{sec: case of r eq 4}\label{sec:fourpoints}

In this section, we study the case of cyclic covers branched at four points. By Lemma \ref{lem: inequality} and Section \ref{sec: reduction}, we have reduced the two main theorems to the following proposition.

\begin{proposition}
\label{prop:main  r=4}
Let $(m, 4, \au)$ be a cyclic monodromy datum, with four branched points, and $p$ a prime, $p\nmid m$.
\begin{itemize}
    \item 
Then, for $p>3m$,  the non-$\D$-ordinary locus of $\M(m,4,\au)_{\fpb}$ is non-empty. 

\item Assume (A2): there exists $\tau'\in\tg$ such that 
\[f(\tau')=1,\qquad f({\tau'}^*)=1, \qquad f(p\tau')=1\]

Then, for $p\geq 3m$,  the number of non-$\D$-ordinary points in  $\M(m, 4,\au)(\fpb)$ is at least $ \lfloor \frac{p}{m}-1\rfloor $. 
\end{itemize}
\end{proposition}

With $\tn$ as in Notation \ref{notation: notation one}, we write $\tau=\tau_0$, and $\OO=\OO_{\tau}=\{p^i\tau\mid 0\leq i<l(\OO)\}$. When (A2) holds, we choose $\tau=\tau_0=\tps$, where $\tp$ satisfying (A2); hence $f(p\tau^*)=1$. With notations as in Notation \ref{notation:t}, we regard the $\D$-ordinary {$\OO$-partial} Hasse--Witt invariant $h_0$ from Definition \ref{def:h_0} as a polynomial $h_0(t)\in \fpb[t]$
and introduce  the $\D$-ordinary {$\tau$-partial} Hasse--Witt invariant $h_1$ as a divisor $h_1(t)$ of $h_0(t)$.

\begin{definition}\label{def:h_1} With the above notations,
     Set $i_0=\min\{1\leq i\leq l \mid f(p^i\ts)=1\}$, where $l = l(\OO)$. We define the $\D$-ordinary {$\tau$-partial} Hasse--Witt invariant $h_1$ as\footnote{This is the $p^{-i_0}\tau$-partial Hasse--Witt invariant in the sense of \cite{Goldring2013TheH}}
\begin{equation}\label{h11}
h_1(t) = A_{i_0-1}(t)\circ\cdots \circ A_0 (t),\end{equation}
where $A_i(t)$ denotes the specialization at $(x_1,x_2,x_3,x_4)=(\infty, t, 1, 0)$ of the matrices $A_i$ in Equation \eqref{eq: Eq_A_i}. 
\end{definition}
By definition,  the polynomial $h_1(t) \in\fpb[t]$ divides 
$ h_0(t)$. 

\begin{remark}\label{partialhasse}
    In \cite[Theorem 7.2]{lin2024abelian}, the $\D$-ordinary Hasse--Witt invariant is constructed as a product of $\tau$-partial Hasse--Witt invariants, as $\tau$ varies in $\tg$. By definition, the $\tau$-partial and $\tau^*$-partial Hasse--Witt invariants agree; furthermore, if $f(\tau)=0$, then the $\tau$-partial Hasse--Witt invariant is trivial.
   By Lemma \ref{lem: ft+fts<=r-2}, when $r\leq 5$, the $\tau$-partial Hasse--Witt invariant is trivial unless $\tau\in \tg$ satisfies $g(\tau)\geq 2$ and either $f(\tau)=1$ or $f(\tau^*)=1$; when $r=4$ that is unless $f(\tau)=1$ and $f(\tau^*)=1$.
   
  Hence, when $r\leq 5$, the $\D$-ordinary Hasse--Witt invariant is a product of $\OO$-partial Hasse--Witt invariants $h_0$ as in Definition \ref{def:h_0}, where $\OO$ varies among the Frobenius orbits in $\tg$ satisfying $\OO=\OO_{\tn}$ for some $\tau_0\in\tg $ satisfying $f(\tau_0)=1$ and $g(\tn)\geq 2$, up to complex conjugation. Furthermore, when $r=4$, each $\D$-ordinary $\OO$-partial Hasse--Witt invariant is a product of $\tau$-partial Hasse--Witt invariants $h_1$ as in Definition \ref{def:h_1}, as $\tau$ varies among the elements in $\OO$ satisfying $f(\tau)=1$ and $f(\tau^*)=1$, up to complex conjugation. 
\end{remark}

By Lemma \ref{lemma:reduction} and Definition \ref{def:h_1}, we have reduced the proof of Proposition \ref{prop:main r=4}  to the following statement. 

\begin{proposition}\label{2stat}
Let $h_1(t)\in \fpb[t]$ the $\D$-ordinary $\tau$-partial Hasse--Witt invariant from Equation \eqref{h11}, and define
 $g_1(t)=t^{-v}h_1(t)$ for $v$ the order of vanishing of $h_1(t)$ at $0$.

Then, the polynomial $g_1(t)$ satisfies $g_1(0)\neq 0,  g_1(1)\neq 0$, and $\deg(g_1)>0$.

    Furthermore, if $f(p\ts)=1$, then  $g_1(t)$ is separable, and $\deg(g_1) \geq  \lfloor \frac{p}{m} -1 \rfloor $.
\end{proposition}

\begin{proof}
Suppose $f(p\ts)=1$ (that is,  $i_0=1$ in Definition \ref{def:h_1}). Then $h_1(t)=A_0(t)=\phi_{\tau}(1,1)(t)$. 

By Corollaries \ref{Cor: deg} and \ref{Cor:vt},  
\[
v_{t}(h_1) 
=\max\big\{0,\st-p+1-\pbto-\pbtd\big\},
\] 
\[
\deg(h_1)= \max\big\{\pbtt , \st-p+1-\pbto \big\}.
\] 
Hence,
$\deg(g_1)=\deg(h_1)-v_t(h_1)$ is equal to one of the following values
\[
\lfloor{p\langle \tfrac{\tau a_2}{m}\rangle}\rfloor,
\lfloor{p\langle \tfrac{\tau a_3}{m}\rangle}\rfloor, 
p-1
-\lfloor{p\langle \tfrac{\tau a_1}{m}\rangle}\rfloor, \text{ or }
p-1
-\lfloor{p\langle \tfrac{\tau a_4}{m}\rangle}\rfloor.  
\]
{By Lemma \ref{lem: inequality}, we have $\frac{p}{m}-1 
\leq \lfloor{p\langle \tfrac{\tau a(k)}{m}\rangle}\rfloor 
\leq p\frac{m-1}{m}$}. Hence,  all of them are greater or equal to $\tfrac{p}{m}-1 $; we deduce $\deg(g_1) \geq\lfloor{\frac{p}{m}-1}\rfloor$.

By Corollary \ref{prop: separability of phit and psit}, the polynomial $g_1(t)$ is separable; furthermore,
$g_1(0)\neq 0$ by definition, and $g_1(1)\neq 0$
since,  by Corollary \ref{Cor: mult of one}, the inequality $\pbtt+\pbtd\leq  p-1$ implies  $h_1(1) \neq 0$.

Suppose $f(p\ts)\neq 1$ (that is, $i_0\geq 2$ in Definition \ref{def:h_1}). Then, to conclude it suffices to show that $h_1(1) \neq 0$ (and hence $g_1(1)\neq 0$) and $\deg(h_1(t))-v_t(h_1(t))>0$ (and hence $\deg g_1(t)> 0$).
We establish these results in Propositions \ref{prop: Better bound when r=4(2): part 2} and \ref{prop: Better bound when r=4(2): part 3}.
\end{proof}

In the following, we maintain the notations of Definition \ref{def:h_1}. We assume $f(p\ts)\neq 1$ (that is, $i_0\geq 2$). 
Then, for $1 \leq i \leq i_0-1$, by definition  $f(p^i\tau^*)\in\{0,2\}$), and we denote \begin{equation}\label{ti}
    \tau_i=p^i\tau \text{ if }f(p^i\ts)=2, \text{ and } \tau_i=p^i\ts  \text{ if }f(p^i\ts)=0. 
\end{equation}


\begin{lemma}\label{aotphi}
Suppose $f(p\tau^*)=2$.
Then 
$A_0(1)
=\phi_{\tau}(1,1)(1) 
\begin{bmatrix} 1 \\ -\frac{1+\pbtfinl}{1+\pbtoinl} \end{bmatrix}$ and 
$\phi_{\tau}(1,1)(1)\neq 0$.
\end{lemma}
\begin{proof}
Since $f(\ts)=1, f(p\ts)=2$, we have $s_\tau=p(f(\tau^*)+1)-(f(p\tau^*)+1)=2p-3$.
By Corollary \ref{Cor: mult of one}, we have $\phi_{\tau}(1,1)(1)\neq 0$. By Lemma~\ref{Lem: value at one}, we see that 
\[
\tfrac{\phi_{\tau}(2,1)(1)}{\phi_{\tau}(1,1)(1)} 
= \tfrac{(-1)^{p-1} \binom{\pbttinl+\pbtdinl}{p-1-\pbtoinl}}{(-1)^{p-2} \binom{\pbttinl+\pbtdinl}{p-2-\pbtoinl} }
\equiv -\tfrac{-1-\pbtf}{-1-\pbto} \pmod{p}. \qedhere 
\]

\end{proof}

\begin{lemma}\label{aotpsi}
Suppose $f(p\tau^*)=0$.
Then $A_0(1)=\psi_{\tau}(1,1)(1)\begin{bmatrix}1 \\ -1 \end{bmatrix}$ and $\psi_{\tau}(1,1)(1)\neq 0$.
\end{lemma}
\begin{proof}
By Corollary \ref{Cor: mult of one}, we know that $\psi_{\tau}(1,1)(1)\neq 0$.
By Lemma \ref{Lem: value at one}, $\psi_{\tau}(2,1)(1)=-\psi_{\tau}(1,1)(1)$. 
\end{proof}

\begin{lemma}\label{aktphi}
Let $1 \leq i \leq i_0-2$.  Suppose $f(p^i\tau^*)= f(p^{i+1}\tau^*)$.
Then 
$A_i(1)=A_i(1,1)(1) \begin{bmatrix} 1 & -1 \\ -b_i & b_i \end{bmatrix}$, where $A_i(1,1)(1)\neq 0$ and
 $ b_i=\bifor $.
\end{lemma}

\begin{proof}
By Equation \eqref{eq: Eq_A_i}, $A_i$ is the matrix of $\phi_{\tau_i}$; hence, by Lemma~\ref{Lem: value at one}, for $1\leq  j, j'\leq 2$,
\[
A_i(j',j)(1)
=(-1)^{\sti-pj+j'}
\binom{\pbtit+\pbtid}{\sti-pj+j'-\pbtio}.
\]
By Corollary \ref{Cor: mult of one}, we have $A_i(1,1)(1)\neq 0$.
By direct computations,
\[\tfrac{A_i(2,1)(1)}{A_i(1,1)(1)} =- \tfrac{(\pbtit+\pbtid)-(\sti-p+1-\pbtio)}{\sti-p+1-\pbtio+1}\equiv -\tfrac{-1-\pbtif}{-1-\pbtio}\pmod{p},\]
\[\tfrac{A_i(2,2)(1)}{A_i(1,2)(1)} =- (\tfrac{\pbtit+\pbtid)-(\sti-2p+1-\pbtio)}{\sti-2p+1-\pbtio+1}\equiv -\tfrac{-1-\pbtif}{-1-\pbtio}\pmod{p},\]
\[
\tfrac{A_i(1,2)(1)}{A_i(1,1)(1)} 
=-\tfrac{\lfloor\frac{\sti-p+1-\pbtio}{p}\rfloor}
{1+\lfloor\frac{\pbtit+\pbtid-(\sti-p+1-\pbtio)}{p}\rfloor}\equiv-\tfrac{1}{1+0}=-1 \pmod{p}.
\]
In the latter equation, the first equality follows from Lemma \ref{choose p less}. 
\end{proof}
\begin{lemma}\label{aktpsi}
Let $1 \leq i \leq i_0-2$. 
Suppose $f(p^i\ts)\neq f(p^{i+1}\ts).$ 
Then 
$A_i(1)=A_i(1,1)(1) \begin{bmatrix} 1 & -1 \\ -1 & 1 \end{bmatrix}$ and $A_i(1,1)(1)\neq 0$.
\end{lemma}
\begin{proof}

Lemma \ref{Lem: value at one} 
implies $A_i(2,1)(1)=-A_i(1,1)(1)$ and $A_i(2,2)(1)=-A_i(1,2)(1)$.
By assumptions, $\sti=3p-3$; following the proof of Lemma \ref{aktphi},
we deduce $A_i(1,1)(1)\neq 0$ and $\frac{A_i(1,2)(1)}{A_i(1,1)(1)} =-1$.
\end{proof}

\begin{lemma}\label{aiphi,aipsi}
$A_{i_0-1}(1)=A_{i_0-1}(1,1)(1)\begin{bmatrix} 1 & -1 \end{bmatrix}$ and  $A_{i_0-1}(1,1)(1)\neq 0$.
\end{lemma}
\begin{proof}
The equality $A_{i}(1,2)(1)=-A_{i}(1,1)(1)$ follows from Lemma \ref{Lem: value at one}; following the proofs of Lemmas \ref{aktphi} and \ref{aktpsi}, with $f(p^{i+1}\tau^*)=1$, we deduce $A_{i_0-1}(1,1)(1)\neq 0$.
\end{proof}

\begin{proposition}\label{prop: Better bound when r=4(2): part 2} Suppose $f(p\ts)\neq 1.$
    Then $h_1(1) \neq 0$.
\end{proposition}
\begin{proof}

By definition,  $h_1(t)= A_{i_0-1}(t)\circ\cdots \circ A_0 (t)$.  By
 Lemmas \ref{aotphi}, \ref{aotpsi}, \ref{aktphi}, \ref{aktpsi}, \ref{aiphi,aipsi},  
 we deduce
\[h_1(1) =\Big(\prod_{i=0}^{i_0-1} A_i(1,1)(1)\Big) \begin{bmatrix} 1 &-1 \end{bmatrix} \begin{bmatrix} 1 &-1\\ -c_{i_0-2} &c_{i_0-2} \end{bmatrix} \dots \begin{bmatrix} 1 &-1\\ -c_{1} &c_{1} \end{bmatrix} \begin{bmatrix}1\\ -c_0\end{bmatrix},  \]
where 
$c_0=\frac{1+\pbtfinl}{1+\pbtoinl}$ if $f(p\ts)=2$, and 
$c_0=1$ if $f(p\ts)=0$, 
and 
$1\leq i\leq i_0-2$, 
\[c_i= \begin{cases}
    \bifor &\text{if } f(p^i\tau^*)=f(p^{i+1}\ts),\\
    1   &\text { if } f(p^i\tau^*)\neq f(p^{i+1}).\end{cases}\]
    Furthemore, $\prod_{i=0}^{i_0-1} A_i(1,1)(1)\neq 0$. 
    By direct computation, 
\[\begin{bmatrix} 1 &-1 \end{bmatrix} \begin{bmatrix} 1 &-1\\ -c_{i_0-2} &c_{i_0-2} \end{bmatrix} \dots \begin{bmatrix} 1 &-1\\ -c_{1} &c_{1} \end{bmatrix} \begin{bmatrix}1\\ -c_0\end{bmatrix} =\prod_{i=0}^{i_0-2} (1+c_i).\]
Hence, it suffices to check $1+c_i\not\equiv 0\pmod{p}$, for each $0\leq i\leq i_0-2$. 
When $c_i=1$, the statement holds trivially; so assume that $c_i=\bifor$.

Let $0\leq b_1,b_4\leq m-1$ satisfying $p\ti a(1)\equiv b_1\bmod{m} \text{ and } p\ti a(4)\equiv b_4\bmod{m}$. 

Then, for $j\in\{1,4\}$,
\[
\lfloor{p\langle \tfrac{\ti a(j)}{m}\rangle}\rfloor 
= p\langle \tfrac{\ti a(j)}{m}\rangle 
-\langle \tfrac{p\ti a(j)}{m}\rangle 
\equiv -\tfrac{b_1}{m} \pmod{p}.
\]
Thus, $1+c_i\not\equiv 0\pmod{p}$ if and only if 
\[
0\not \equiv 
2+\pbtio+\pbtif
\equiv\tfrac{2m-b_1-b_4}{m} \pmod{p}.
\]
Since $2\leq 2m-b_1-b_4\leq 2m<p$, we deduce $ 2m-b_1-b_4\not\equiv 0\pmod{p}$. 
\end{proof}

\newcommand{\fR}{{\mathfrak R}}

Similarly to the computation of the $\D$-ordinary partial Hasse--Witt invariant $h_0$ in Equation \eqref{RJi}, following \cite[Equations (19) and (20)]{lin2024abelian}, we compute the polynomial $h_1$ in Equation \eqref{h11} as 

\begin{equation}\label{RJi1}
h_1(t)=\sum_{J\in \jj_0} \fR_J(t),
\text{ with }
\fR_J(t)=\prod_{i=0}^{i_0-1} \fR_{J,i}^{p^{i_0-1-i}}.  
\end{equation}
where $\mathfrak{J}_0$ denote the set $$\mathfrak{J}_0=\{J:\{0,1, \dots, i_0\} \to \mathbb{N} \mid 
J(0)=J(i_0)=1,  1 \leq J(i) \leq d(i) \text{ for $1 \leq i \leq i_0$}\}$$
with $d(i)$ are defined as in \eqref{def: ai}, and
for any $J \in \mathfrak{J}_0$, ${\mathfrak R}_{J,i}=A_i(J(i+1),J(i))$.


\begin{lemma}\label{L1}
For any $J_1,J_2\in\jj_0$: 
\begin{enumerate}
    \item if $v_t(\fR_{J_1}(t))=v_t(\fR_{J_2}(t))$, then $J_1=J_2$.
    \item if $\deg(\fR_{J_1}(t))=\deg(\fR_{J_2}(t))$, then $J_1=J_2$.
\end{enumerate}
\end{lemma}
\begin{proof}
 We prove (1) (the proof of (2) is similar). By definition,
\[v_t(\fR_{J_k}(t))=\sum_{i=0}^{i_0-1}p^{i_0-1-i}v_t(A_{i}(J_k(i+1),J_k(i))(t)),\]
By Lemmas \ref{aotphi}, \ref{aotpsi}, \ref{aktphi}, \ref{aktpsi}, \ref{aiphi,aipsi}, we have $0\leq v_t(A_{i}(J_k(i+1),J_k(i))(t))<p$,  for all $0\leq i\leq i_0-1$.
Hence, $v_t(\fR_{J_1}(t))=v_t(\fR_{J_2}(t)) $ implies $$ v_t(A_{i}(J_1(i+1),J_1(i))(t))=v_t(A_{i}(J_2(i+1),J_2(i))(t))\text { for all } 0\leq i\leq i_0-1.$$ 

Now we prove $J_1(i)=J_2(i)$ by reverse induction on $i$. By definition, $J_1(i_0)=J_2(i_0)=1$. Let $1\leq i\leq i_0-1$, and assume $J_1(i+1)=J_2(i+1)$. 
Let $\tau_i$ as in Equation \eqref{ti},
hence $f(\tau_i^*)=2$. By Lemma \ref{Cor:vt}, if  $J_1(i+1)=J_2(i+1)$, the equality $v_t(A_{i}(J_1(i+1),J_1(i))(t))=v_t(A_{i}(J_2(i+1),J_2(i))(t))$ implies $J_1(i)=J_2(i)$.
\end{proof}

\begin{proposition}\label{prop: Better bound when r=4(2): part 3} Suppose $f(p\ts)\neq 1$.
   Then, $\deg(h_1(t))-v_t(h_1(t))>0$.
\end{proposition}

\begin{proof} By Equation \eqref{RJi1}, we have
$v_t(h_1(t))\geq \min_{J\in\jj_0} v_t(\fR_J(t)),$ and 
$\deg(h_1(t))\leq \max_{J\in\jj_0} \deg(\fR_J(t)).$

By Lemmas \ref{L1}, there exist a unique $J_1\in\jj_0$ that minimizes $v_t(\fR_{J_1}(t))$, and a unique $J_2\in\jj_0$ that maximizes $\deg(\fR_{J_2}(t))$). We deduce
\[v_t(h_1(t))=v_t(\fR_{J_1}(t))< \deg(\fR_{J_1}(t))\leq \deg(\fR_{J_2}(t))=\deg(h_1(t)).\]
The  inequality $v_t(\fR_{J_1}(t))< \deg(\fR_{J_1}(t))$ follows from Corollary \ref{Cor: bound on degree} since $p>3m$.
\end{proof}

\section{Non-$\D$-ordinary covers branched at more than four points}\label{sec: main case} 

In this section, we study the case of cyclic covers branched at more than four points. By Section \ref{sec: reduction}, we have reduced the Theorem  \ref{main theorem} to the following proposition.
\begin{proposition}\label{main prop}

Let $(m,r,\au)$ be a cyclic monodromy datum, with $r\geq 4$. Suppose there exists $\tn \in \tg$ such that 
$f(\tn)=1$ and $g(\tn)\geq 2$.

Then,  for any prime $p\geq (3\binom{r}{2}+r-2)m$,  the non-$\D$-ordinary locus of $\M(m,r,\au)_{\fpb}$ is non-empty.

\end{proposition}

\begin{proof}
We follow the same strategy as in the proof of Proposition \ref{prop:main  r=4}. Let $\D$-ordinary $\OO_{\tau_0}$-partial Hasse--Witt invariant $h_0\in \fpb[x_1, \dots, x_r]$ defined in Equation \eqref{h_0}. 
By Lemma \ref{lemma:reduction}, it suffices to show that, for each pair $1\leq i<j\leq r$, the multiplicity of the boundary divisor $x_i-x_j$ in $h_0$ is strictly less than the degree of $h_0$. By Lemma \ref{degree}, it suffices to show 
\begin{equation}\label{eqpropmain}
    \sum_{1\leq j_1<j_2\leq r}v_{x_{j_1}-x_{j_2}}(h_0)< p^{l-1}(p-(r-2)),
    \end{equation}
where $v_{x_{j_1}-x_{j_2}}(h_0)$ denotes the multiplicity of $x_{j_1}-x_{j_2}$ as divisor of $h_0$. 
We establish this inequality in Proposition \ref{main proposition both}.
\end{proof}

Note that in Proposition \ref{main prop}, the bound on $p$ is greater than that in Proposition \ref{prop:main  r=4}, as the latter was  obtained by computing --and not just bounding-- the multiplicities of the boundary divisors and the degree of (a divisor of) $h_0$.

Let $f(y)\in R[y]$ be a single variable polynomial in $y$ with coefficients in a ring $R$.
If $R$ has characteristic $0$, then one can compute $v_y(f)$ by studying the derivatives of $f$, that is  $v_y(f)=t$ if and only if $f^{(t)}(0) \neq 0$ and $f^{(t')}(0)$ for all $0\leq t'<t$. If $R$ has positive characteristic, the analogous statement is false. To bound $v_y(f)$ when $R$ is a ring of characteristic $p$, we introduce a suitable family of differential operators.

\begin{definition}\label{Dpk}
Let $R$ be a polynomial ring over a field of characteristic $p$. For any $k \geq 0$,  define $D_{p^k}:R[y]\to R[y]$ as follows, for
$f(y)=\sum a_iy^i\in R[y]$,
\begin{align}\label{eqdpk}
D_{p^k}(f(y)) \coloneqq 
\sum_{i} \lfloor \tfrac{i}{p^k}\rfloor 
a_iy^{i-p^k}.
\end{align}
\end{definition}

By definition, the operators $D_{p^k}$ are $R$-linear and $D_1(f(y))=\frac{d}{dy}f(y)$.

\begin{lemma}\label{property of Dpk1}
For any $k\geq 0$, and $f(y)\in R[y]$, we have $$   v_y(f(y))  \leq p^k+v_{y}(D_{p^k}(f(y))). $$

   
\end{lemma}
\begin{proof}
Let $f(y)=\sum a_iy^i\in R[y]$, and denote $i_0= v_{y}(f(y))$ and $i_1= v_{y}(D_{p^k}(f(y))$. By definition, 
$i_0=\min\{i\mid a_i \neq 0\}$ and $i_1=\min\{i\mid  \lfloor\frac{i+p^k}{p^k}\rfloor  a_{i+p^k} \neq 0 \}$. In particular,  $a_{i_1+p^k}\neq 0$, and hence $i_0\leq i_1+p^k$. 
\end{proof}


\begin{lemma}\label{property of Dpk2}
For any $k'> k\geq 0$, and any  $f_1, f_2, f_3 \in R[y]$ satisfying $\deg_x(f_3)< p^{k}$, we have
$$D_{p^k}(f_1^{p^{k'}}f_2^{p^k}f_3)=f_1^{p^{k'}}(f_2')^{p^k}f_3. $$
\end{lemma}
\begin{proof}
Let $f_1=\sum a_i y^i$, $f_2=\sum b_j y^j$, $f_3=\sum c_s y^s$, and hence
$f_1^{p^{k'}}f_2^{p^k}f_3=\sum_{i,j,s}a_i^{p^{k'}}b_j^{p^k}c_s y^{p^{k'}i+p^kj+s}$. 

For any $i, j, s\geq 0$, if  $s<p^k$, then
$\lfloor\frac{p^{k'}i+p^kj+s}{p^k}\rfloor \equiv j\bmod p$.
Hence, by assumptions,
 $$D_{p^k}(f_1^{p^{k'}}f_2^{p^k}f_3) 
    =\sum_{i,j,s}a_i^{p^{k'}}b_j^{p^k}c_s y^{p^{k'}i+p^k(j-1)+s} 
    =\sum_{i,j,s} (a_iy^{i})^{p^{k'}} (jb_jy^{j-1})^{p^k} c_s y^s =f_1^{p^{k'}}(f_2')^{p^k}f_3.
    $$
    \qedhere
\end{proof}


In the following, 
we apply the differential operators from Equation \eqref{eqdpk} to bound the multiplicty of the boundary in the $\D$-ordinary partial Hasse--Witt invariant.


\begin{notation}\label{notation: notation 2}
    We fix a cyclic monodromy datum $(m,r,\au)$, and an element $\tau \in \tg$ such that $f(\ts)=1$ and $g(\ts)\geq 2$, and a prime $p$ such that $p\geq (3\binom{r}{2}+r-2)m$. 
    
    We denote $\OO=\OO_\tau$ the Frobenius orbit of $\tau$, of size $l=l(\OO_\tau)$, and fix the bijection $\{0,1, \dots, l-1\} \to \OO$ via $i \mapsto p^i\tau$.  For $1 \leq k \leq r$, we set $b_k=\lfloor p \langle\frac{\tau a(k)}{m}\rangle \rfloor$. 

We write $h_0$ for the $\D$-ordinary $\OO$-partial Hasse--Witt invariant from Equation \eqref{h_0}.
\end{notation}

Given $1\leq j_1<j_2\leq r$, to study $v_{x_{j_1}-x_{j_2}}(h_0)$, we set  $y=x_{j_1}-x_{j_2}$ and, after a change of variable, regard a polynomial $f(x_1, \dots x_r)$ as  a polynomial $f(y)=f(x_1, \dots, \hat{x}_{j_2}, \dots, x_r)(y)$ of a single variable $y$ with coefficients in the ring $R=\fp[x_1, \dots, \hat{x}_{j_2}, \dots, x_r]$, and write  $f^{(t)}=\frac{\partial^t}{\partial y^t}f(y)$, for $t\geq 0$.

\begin{proposition}\label{main proposition both}
Let $h_0$ as in Equation \eqref{h_0}. Then, for any $1\leq j_1<j_2\leq r$, we have
$v_{x_{j_2}-x_{j_1}}(h_0) \leq \sum_{i=0}^{l-1}p^{l-1-i} \cdot t_i'$
where for $1\leq i\leq l$, 
$$
t_i'= \begin{cases}
\max\{0,\lfloor p \langle \frac{p^i\tau a(j_1)  }{m}\rangle\rfloor+\lfloor p \langle \frac{p^i\tau a(j_2)}{m}\rangle \rfloor-(p-3)\} &\text{if } f(p^i\ts)\geq 1,\\
\max\{0, \lfloor{p\langle\frac{-p^i\tau a(j_1)}{m}\rangle}\rfloor+
 \lfloor{p\langle\frac{-p^i\tau a(j_2)}{m}\rangle}\rfloor -(p-3)\} &\text{if } f(p^i\ts)=0.
 \end{cases}
$$
In particular, $\sum_{1\leq j_1<j_2\leq r} v_{x_{j_2}-x_{j_1}}(h_0) <p^{l-1}(p-(r-2)).$
\end{proposition}

\begin{proof}
Without loss of generality, set $j_2=r$ and $j_1=r-1$. Let $y=x_{r}-x_{r-1}$, and for $R=\fp[x_1,\dots,x_{r-1}]$ consider the operators $D_{p^k}$ given in Definition \ref{Dpk}  .


Set
\[
t_i
=\begin{cases}
\max\{0,\pbtr+\pbtrm-(p-1)\} &\text{if }   f(p^i\tn^*)\geq 1, f(p^{i+1}\tn^*)\geq 1 ,\\
\max\{0,{\lfloor{p\langle\frac{-p^i\tau a(r)}{m}\rangle}\rfloor}+
\lfloor{p\langle\frac{-p^i\tau a(r-1)}{m}\rangle}\rfloor
-(p-1)\} 
&\text{if } 
f(p^i\ts)=0, f(p^{i+1}\ts)=0;\\
\max\{0,\pbtr+\pbtrm-(p-3)\} &\text{if } 
  f(p^i\ts)\geq 1, f(p^{i+1}\ts)=0 \\
 \max\{0,\lfloor{p\langle\frac{-p^i\tau a(r)}{m}\rangle}\rfloor+
 \lfloor{p\langle\frac{-p^i\tau a(r-1)}{m}\rangle}\rfloor-(p-3)\} 
 &\text{if }  
   f(p^i\ts)=0, f(p^{i+1}\ts)\geq 1.
\end{cases}\]
and
 $   h_2(y) = D_{p^{l-1}}^{t_0} \circ D_{p^{l-2}}^{t_1} \dots D_{p^1}^{t_{l-2}} \circ D_{p^0}^{t_{l-1}}(h_0(y))\in R[y]. $

By definition, $t_i\leq t_i'$;
hence, by Lemma \ref{property of Dpk1}, we have
\[v_{x_r-x_{r-1}}(h_0)=v_y(h_0(y)) \leq v_y(h_2)+\sum_{i=0}^{l-1}p^{l-1-i}t_i \leq v_y(h_2)+\sum_{i=0}^{l-1}p^{l-1-i}t_i'.\]

Hence, to deduce the stated bound on $v_{x_r-x_{r-1}}(h_0)$, it suffices to show $v_y(h_2(y))=0$, or equivalently $h_2(0)\neq 0$  in $R$. 

From Equation \eqref{RJi}, we deduce $$ h_2(y)=\sum_{J\in\mathfrak J} D_{p^{l-1}}^{t_0} \circ D_{p^{l-2}}^{t_1} \dots D_{p^1}^{t_{l-2}} \circ D_{p^0}^{t_{l-1}}(\prod_{i=0}^{l-1}R_{J,i}(y)^{p^{l-i-1}})); $$
by definition of $R_{J,i}$, for all $J, i$,  Equations \eqref{coefficient of phii} and \eqref{coefficient of psii} imply $\deg_y(R_{J,i}(y))<p$. Hence,  by Lemma \ref{property of Dpk2}, we deduce
\begin{equation}\label{h2RJi}
    h_2(y)
    =\sum_{J\in\mathfrak J} \prod_{i=0}^{l-1}
    (R_{J,i}^{(t_i)})^{p^{l-i-1}}.
\end{equation}

For any $J\in \mathfrak J$ and $0\leq i<l$, set $M_{J,i} = 
R_{J,i}^{(t_i)}(0)\in R$, and  $M_J=\prod_{i=0}^{l-1}M_{J,i}^{p^{l-i-1}}$. By Equation \eqref{h2RJi}, we have
$h_2(0)=\sum_{J\in \mathfrak{J}} M_{J}.$

 %

Let ${\mathfrak J}_1=\{J\in\mathfrak J\mid M_J\neq 0\}$. By Lemma 
\ref{lem: entry of phi and psi},
if $J\in \mathfrak{J}$ satisfies $J(i)=1$ for all $0\leq i\leq l$, then $M_{J}\neq 0$; hence, $\mathfrak J_1$ is non-empty, and  $h_2(0)=\sum_{J\in {\mathfrak J}_1} M_J$. 

For any  $J\in \mathfrak J_1$ and $0\leq i<l$, denote by $T_{J,i}$ the maximal monomial in $M_{J,i}$, under the lexicographical order on $R$ (given by $x_i>x_j$ for $i<j$), and set  $T_J =\prod_{i=0}^{l-1}T_{J,i}^{p^{l-i-1}}$. Then, $T_J$ is the maximal monomial of $M_J$, and 
to prove $h_2(0)\neq 0$, it suffices to prove that for any $J_1,J_2\in \mathfrak{J}_1$, if $T_{J_1}=T_{J_2}$ then $J_1=J_2$. 

  We prove the latter statement  by adapting a similar argument in the proof of \cite[Theorem 7.2]{lin2024abelian}. 
Let $J_1,J_2\in \mathfrak{J}_1$ satisfying $T_{J_1}=T_{J_2}$.  
Then, 
\[v_{x_s}(T_{J_1})-v_{x_s}(T_{J_2})=\sum_{i=0}^{l-1}p^{l-1-i}(v_{x_s}(T_{J_1,i})-v_{x_s}(T_{J_2,i}))=0.\]

For $1\leq s\leq r-1$ and $0\leq i\leq l-1$, set $\eta_{i,s}=v_{x_s}(T_{J_1,i})-v_{x_s}(T_{J_2,i})$. 
By Lemma 
\ref{lem: entry of phi and psi}, we deduce $|\eta_{i,s}|<p$, and hence  $\eta_{i,s}=0$, for all $1\leq s\leq r-1$ and $0\leq i\leq l-1$. Furthermore,
\[\deg(T_{J_1,i})-\deg(T_{J_2,i}) =\sum_{s=1}^{r-1} \eta_{i,s}=0.\]
By Lemma \ref{lem: entry of phi and psi}, for all $J\in {\mathfrak J}_1$, $\deg(T_{J,i})=
s_{\tau}-pJ(i)+J(i+1)$, 
Hence, by induction on $l-i$, we conclude (by definition, $J(l)=1$ for all $J\in {\mathfrak J}$).
To conclude, we prove
$\sum_{1\leq j_1<j_2\leq r} v_{x_{j_2}-x_{j_1}}(h_0) <\deg(h_0).$

Recall from Notation \ref{notation: notation 2}, $b_k=\ptnk$, for $1\leq k\leq r$, 
and set $s=\sum_{k=1}^r b_k$. 
For any $1\leq j_1<j_2\leq r$, the established bound on $v_{x_{j_1}-x_{j_2}}(h_0)$ implies 
\begin{align*}
v_{x_{j_1}-x_{j_2}}(h_0) 
& \leq p^{l-1} \max\{0, b_{j_1}+b_{j_2}-p+3\}+\tsum_{i=1}^{l-1}p^{l-1-i}(p-1) \\
& < p^{l-1}(2+\max\{0,b_{j_1}+b_{j_2}-p+1\}) +p^{l-1};
\end{align*}
we deduce
\[
\sum_{1\leq j_1<j_2\leq r}v_{x_{j_1}-x_{j_2}}(h_0) 
< p^{l-1}\Big(3 \tbinom{r}{2}+\sum_{j_1<j_2}\max\{0,b_{j_1}+b_{j_2}-p+1\}\Big).
\]
By Lemma \ref{degree}, $\deg(h_0)\geq p^{l-1}(p-(r-2))$; 
by Lemma \ref{inequality}, we
conclude. 
\end{proof}

\begin{lemma}\label{inequality}
With notations as in Notation \ref{notation: notation 2}, where $p>(3\tbinom{r}{2}+r-2)m$, we have
\begin{equation}\label{eq: inequality}
3\tbinom{r}{2}+\sum_{j_1<j_2}\max\{0,b_{j_1}+b_{j_2}-p+1\} \leq p-(r-2))
\end{equation}
\end{lemma}
\begin{proof}
Without loss of generality, we reorder the indices $1\leq k\leq r$ so that $b_1 \geq b_2 \geq \dots \geq b_r$. 
Let $S=\{(j_1,j_2) \mid 1\leq j_1 < j_2\leq r\text{ and } b_{j_1}+b_{j_2} \geq p\}$; for $1\leq k\leq r$, set $S_k=\{j \mid (k,j) \in S\}$ and 
$n=|\{k\mid |S_k|\neq 0\}|$. Then
\[S=\{(k,j)\mid 1\leq k\leq n \text{ and } k+1\leq j\leq |S_k|\}.\]

We claim that $n \leq 2$. Indeed, suppose $n\geq 3$. Then $(3,4)\in S$, that is $b_3+b_4\geq p$, and hence also $b_1+b_2\geq p$. we deduce the contradiction
\[2p\leq b_1+b_2+b_3+b_4 \leq s<2p,\]
where the latter inequality follows from $s=\sum_{k=1}^r b_k=2p-f(p\tau_0^*)-1<2p$.

Assume $n=0$, that is $S=\emptyset$. Then the inequality in \eqref{eq: inequality} follows from $p\geq (3\binom{r}{2}+r-2)m> 3\binom{r}{2}+r-2$.

Assume $1\leq n\leq 2$. Let $V=\sum_{k<j}\max \{0, b_k+b_j-p+1\}$. With the above notations,
\begin{align*}
    V
    &= \sum_{k=1}^n \sum_{j=k+1}^{k+|S_k|}(b_k+b_j-(p-1)) 
    \leq-(p-1)(\sum_{k=1}^n |S_k|)+ \sum_{k=1}^n |S_k|b_k+\sum_{k=1}^n(s-\tsum_{j=1}^k b_j)\\
    &= -(p-1)(\tsum_{k=1}^n |S_k|)+\tsum_{k=1}^n (|S_k|-n-1+k)b_k+n s.
\end{align*}
For each $1\leq k\leq n$, we have $|S_k| \geq n+1-k$, $b_k\geq p\frac{m-1}{m}$ and $s\geq 2p-2$.
Therefore,
\begin{align*}
V& <  -(p-1)(\tsum_{k=1}^n |S_k|)+\tsum_{k=1}^n (|S_k|-n-1+k)(p-\tfrac{p}{m})+2np-2n \\
    &=p\big[ -(n+1)n+\tfrac{n(n+1)}{2}+2n\big]+(\tsum_{k=1}^n |S_k|)-\tfrac{p}{m}\tsum_{k=1}^n (|S_k|-(n+1-k)) -2n, 
\end{align*}
where $-(n+1)n+\frac{n(n+1)}{2}+2n=1$, for $1\leq n\leq 2$. 
We have reduced the proof of inequality \eqref{eq: inequality} to showing 
\[3\tbinom{r}{2}+(r-2)-2n+\tsum_{k=1}^n |S_k|\leq \tfrac{p}{m}\tsum_{k=1}^n (|S_k|-n-1+k).\]
Note that $ |S_k| \geq n+1-k$, for any $k$. 
In particular, $\sum_{k=1}^n (|S_k|-n-1+k) \geq 0$. 

Suppose $n=1$ and $\sum_{k=1}^n (|S_k|-n-1+k)=0$. Then, $|S_1|=1$ and $|S_k|=0$, for any $k>1$. We have
    \begin{align*}
    3\tbinom{r}{2}+V 
    &=3\tbinom{r}{2}+(b_1+b_2-p+1)\leq 3\tbinom{r}{2} +2p(1-\frac{1}{m})-p+1\\
    &\leq 3\tbinom{r}{2}+1+p -2\big(3\tbinom{r}{2}+r-2\big) <p-(r-2).
    \end{align*} 

Suppose $n=2$ and  $\sum_{k=1}^n (|S_k|-n-1+k)=0$. Then $|S_1|=2, |S_2|=1$, and $|S_k|=0$, for any $k>2$. We have
    \begin{align*}
    3\tbinom{r}{2}+V
    &=3\tbinom{r}{2}+2b_1+2b_2+2b_3-3(p-1)\leq 3\tbinom{r}{2}-3(p-1)+2(s_{\tau_0}-b_4)\\
    &\leq 3\tbinom{r}{2}-3(p-1)+2(2p-2-\tfrac{p}{m})\leq 3\tbinom{r}{2} +p-1- 2(3\tbinom{r}{2}+r+2) <p-(r-2). 
    \end{align*}

Suppose $\sum_{k=1}^n (|S_k|-n-1+k)\geq 1$.  Equivalently, $\sum_{k=1}^n |S_k|\geq \frac{n(n+1)}{2}+1=2n$, for $1\leq n\leq 2$.
Then
    \begin{align*}
    &\tfrac{p}{m}\tsum_{k=1}^n (|S_k|-n-1+k) -\tsum_{k=1}^n |S_k|
    = (\tfrac{p}{m}-1)\tsum_{k=1}^n |S_k| -\tfrac{p}{m} \tfrac{n(n+1)}{2}\\
    &\geq (\tfrac{p}{m}-1)(\tfrac{n(n+1)}{2}+1) -\tfrac{p}{m} \tfrac{n(n+1)}{2}=\tfrac{p}{m} -\tfrac{n(n+1)}{2}-1\\
    &\geq \big(3\tbinom{r}{2}+r-2\big) -2n. \qedhere
    \end{align*}
\end{proof}


We conclude this section with a proof of Corollary \ref{cor: cases of main theorem} which identify monodromy data that satisfy Assumption (A1) in Theorem \ref{main theorem}

\begin{lemma}
    \label{Prop: sig 0 then 1}

Let $(G,r,\au)$ be an abelian monodromy datum, of signature type $\fb$. 
Assume there exists $\tau \in \tg$ such that $f(\tau)=0$ and $g(\tau) \geq 2$. 

Then, there exists $\tau_0\in \tg$ such that $f(\tau_0)=1$ and $g(\tau_0) \geq 2$.
\end{lemma}
\begin{proof}
Let $Q=G/\ker(\tau)$.Then, $Q$ is cyclic; we write $m=|Q|$, and denote by $(m, r',\ub)$ the monodromy datum of the quotient $Q$-cover  $C/\ker(\tau)$, of $\Po$.  As  in Section \ref{tg}, we identify $\Z/m\Z =\mathcal{T}_Q$.

By definition, the character $\tau\in \tg$ descends to a character $\bar{\tau}$ in $\mathcal{T}_Q$, satisfying $f_Q(\bar{\tau})=f(\tau)$.  and $g_{Q}(\bar{\tau})=g(\tau)$. Hence, $f_Q(\bar{\tau})=0$  and $g_{Q}(\bar{\tau})\geq 2$.  
Denote by $A\in\Z$, $0\leq A<m$,  the integer corresponding to $\bar{\tau}^*$ under the identification $\Z/m\Z =\mathcal{T}_Q$. By Equation \eqref{signature},  $f_Q(\bar{\tau})=0$ implies $\sum_{k=1}^{r'} \langle\tfrac{Ab(k)}{m} \rangle=1.$

Consider the sequence $f(\bar{\tau}), f(\bar{\tau}^{2}),\dots, f(\bar{\tau}^{m-1})$. By assumption, $f(\bar{\tau}^{m-1})=f(\bar{\tau}^*)=g(\bar{\tau})\geq 2$. Moreover, for any $1\leq t\leq m-2$,
\begin{align}\label{ft}
f(\bar{\tau}^{t+1}) -f(\bar{\tau}^{t}) = \tsum_{k=1}^{r'} \langle\tfrac{(t+1)Ab(k)}{m} \rangle -\langle\tfrac{tAb(k)}{m} \rangle \leq \tsum_{k=1}^{r'} \langle\tfrac{Ab(k)}{m} \rangle =1.
\end{align}
We conclude that $f(\bar{\tau}^{t})=1$  for some $1\leq t\leq m-2$.

Let $t_1=\min\{1\leq t\leq m-2\mid f(\bar{\tau}^{t})=1 \}$. By definition,  $t_1\geq 2$ and $f(\bar{\tau}^{t_1-1})=0$. By Equation \eqref{ft},  $\langle\tfrac{t_1 A b(k)}{m} \rangle =\langle\tfrac{(t_1-1)Ab(k)}{m} \rangle +  \langle\tfrac{Ab(k)}{m} \rangle$, for all $1\leq k\leq r'$.
In particular, if $\langle\frac{t_1 A b(k)}{m} \rangle=0$  then $\langle\frac{A b(k)}{m} \rangle=0$, which implies $g(\tau^{t_1})\geq 2$ since
\[g(\tau^{t_1}) =r'-2-|\{1\leq k\leq r'\mid \langle\tfrac{t_1 A b(k)}{m} \rangle=0\} | \geq r'-2-|\{1\leq k\leq r'\mid \langle\tfrac{A b(k)}{m}  \rangle =0\}| = g({\tau}).\qedhere\]
\end{proof}

By Proposition \ref{Prop: sig 0 then 1}, we identify many instances when Theorem \ref{main theorem} applies.

\begin{proof}[Proof of Corollary \ref{cor: cases of main theorem}] We verify  Assumption (A1) of Theorem \ref{main theorem} by Lemma \ref{Prop: sig 0 then 1}.
\begin{enumerate}
    
\item Assume $4\leq r\leq 5$: for any $\tau \in \tg$, $f(\tau)+f(\ts) \leq r-2 \leq 3$, hence, either $f(\tau)\leq 1$ or $f(\ts)\leq 1$. 
Furthermore, if  $g(\tau) \leq 1$ for all $\tau \in \tg$, then  $\dim \shgf=0$ which contradicts $ \dim(\mgra)=r-3 \geq 1$. So $g(\tau)\geq$ for some $\tau\in \tg$.

\item  Let $m=|G|$ , and $\tau\in tg $ denotes the character corresponding to $1\in\Z$, under the identification $\tg=\Z/m\Z$ from Section \ref{tg}. Then, $g(\tau)=g(\tau^*)=r-2\geq 2$. Denote $S=\sum_{i=1}^r a(i)$. By Equation \ref{signature}, 
$f(\tau)=0$ if $S=m$,  $f(\tau)=1$ if $S=2m$, $f(\tau^*)=0$ if $S=(r-2)m$, and  $f(\tau^*)=1$ if $S=(r-1)m$.

\item  Let $H$ be the unique subgroup $G$ of index $m_1$. For any   $\tau\in \tg$ satisfying $\ker(\tau)=H$,
\[
g(\tau)
=-2+\sum_{i=1}^r\langle\tfrac{a(i)}{m_1}\rangle
+\sum_{i=1}^r\langle\tfrac{-a(i)}{m_1}\rangle
=-2+|\{i \mid m_1 \nmid a(i)\}| \in \{2,3\}.
\] 
Hence, $g(\tau) =f(\tau)+f(\ts)\geq 2$ and either $f(\tau)=1$ or $f(\ts)=1$ (or both). \qedhere
\end{enumerate}
\end{proof}
\begin{remark}
    By \cite[Section 6]{li2019newton},  Corollary \ref{cor: cases of main theorem} gives rise to inductive systems of monodromy data, with unbounded numbers of branched points (and not satisfying assumption (A1) in Theorem \ref{main theorem}), for which there exist smooth non-$\D$-ordinary covers of $\Po$. More precisely, given a pair of cyclic monodromy data $(\gamma_1,\gamma_2)$, of the same degree and branched at at most five points, that is admissible and balanced at $p$ in the sense of \cite{li2019newton}, by \cite[Theorem 5.4]{li2019newton} and \cite[Theorem 1.1]{lin2024abelian}, Theorem \ref{prop: main theorem for r eq 4} impies 
     the existence a non-$\D$-ordinary cover of $\Po$ with monodromy datum equal to that obtained by clutching $(\gamma_1, \gamma_2)$ (see \cite[Definition 3.5]{li2019newton}).  
     Furthermore, by Proposition \ref{prop: other fam}(2), for $m=5$ or $m=7$, (some of) these inductive systems consists of monodromy data for which there exist smooth almost-$\D$-ordinary covers of $\Po$. 
We refer to \cite[Remarks 8.5, 8.6]{lin2024abelian} for two examples.

\end{remark}

\section{An example: the family $\M(7,4,(3,1,1,2))$}\label{sec: example of a family}

In this section, we illustrate by an example how, given an explicit monodromy datum with four branched points,  
Theorem \ref{prop: main theorem for r eq 4} can be strengthened to identify non-empty Ekedahl-Oort and Newton strata in the mod $p$ reduction of the corresponding Hurwitz space.

In the following, we call a cover of $\Po$ {\bf almost-$\D$-ordinary} (resp. {\bf $\D$-special}) if  its \EO type corresponds to a stratum of $\Sh(\D)_{\fpb}$ of codimesion 1 (resp. dimension 0), that is to an almost-$p$-ordinary (resp. $p$-special) \EO stratum.

\begin{lemma}(\cite[Theomre 7.2]{lin2024abelian})\label{lem:almost}
Let $(m, 4, \au)$ be a cyclic monodromy datum, with four branched points, and $p$ a prime, $p\nmid m$. For each $\tau\in \tg$ satisfying $f(\tau)=1$ and $f(\tau^*)=1$, denote by $h_\tau$ the associated $\D$-ordinary $\tau$-partial Hasse--Witt invariant $h_1$ from Definition \ref{def:h_1}. 

If  $\alpha\in \fpb$, $\alpha\neq 0, 1$,  satisfies $h_{\tau_0}(\alpha)=0$ for some $\tau_0\in\tg$ satisfying $f(\tn)=1$ and $f(\tns)=1$, and $h_\tau(\alpha)\neq 0$ for all $\tau\in\tg$, $\tau\neq \tau_0,\tau_0^*$, satisfying $f(\tau)=1$ and $f(\ts)=1$, then 
$\alpha$ defines an almost-$\D$-ordinary point in $\M(m,r,\au)(\fpb)$.   

If $\alpha\in \fpb$, $\alpha\neq 0$, satisfies $h_{\tau}(\alpha)=0$, for all $\tau\in \tg$, then 
$\alpha$ defines a $\D$-special point in $\M(m,r,\au)(\fpb)$. 
\end{lemma}
\begin{proof}
    
As discussed in Remark \ref{partialhasse}, when $r=4$,  by \cite[Theomre 7.2]{lin2024abelian}, the $\D$-ordinary Hasse--Witt invariant is constructed as a product of  
the $\tau$-partial Hasse--Witt invariants $h_\tau$, as $\tau$ varies among the characters satisfying $f(\tau)=1$ and $f(\ts)=1$, up to complex conjugation.  By construction, for each  $\tn$,  a point in $V_{\tn}\setminus \bigcup_{\tau\neq \tn,\tns} V_\tau\cap V_{\tn}$, (if it exists) corresponds to an almost-$\D$-ordinary cover, and a point in $\bigcap_{\tau\in\tg} V_\tau$  to a $\D$-special point (here, $V_\tau=V(h_\tau)$ denotes the vanishing locus of $h_\tau$).
\end{proof}

We carry out such an analysis in the case of the monodromy datum $(7,4,(3,1,1,2))$; the results are summarized in Proposition \ref{prop: special case}.

\begin{notation}
Fix the monodromy datum $\gamma=(7,4,(3,1,1,2))$. Denote $\M=\M(\gamma)$, $\D=\D(\gamma)$ and $\Sh=\Sh(\D)$; let  $p$ be a prime, $p\neq 7$.
\end{notation}

By Equation \eqref{signature type of monodromy},  the signature type  associated with the monodromy $\gamma$ is $\fb=(0,1,1,1,1,2)$; in particular, $\dim\Sh=2$. 

By  \cite[Theorem 1.6]{viehmann2013ekedahl}), we compute the  \EO and Newton strata in the reduction modulo $p$ of $\Sh$.  By \cite[Theorem B]{gortz2019fully}, we deduce that in this instance the  Ekedahl-Oort stratification is a refinement of the Newton stratification.

\begin{lemma}\label{lem:strata}
    Let $p$ be a prime, $p\neq 7$. Then \begin{enumerate}
        \item 
  The non-empty Ekedahl--Oort strata in $\Sh_{\fpb}$ are the open ($p$-ordinary) stratum, two locally closed strata of codimension 1 (almost-$p$-ordinary), and the closed ($p$-special)
    stratum  of dimension 0.


\item The non-empty Newton strata in $\Sh_{\fpb}$ are:
    \begin{enumerate}
    \item if $p \equiv 2, 3, 4, 5 \mod 7$: 
   the open ($\mu$-ordinary) stratum  and the closed (basic) stratum, of codimension 1; 

    \item if $p \equiv 1, 6 \mod 7$: the open  ($\mu$-ordinary) stratum , one locally closed stratum of codimension 1, and the closed  (basic) stratum, of dimension 0.

\end{enumerate}

\item Each \EO stratum is fully contained in a Newton stratum. More precisely,  
the open strata agree, and
\begin{enumerate}
    \item 
 if $p \equiv 2, 3, 4, 5 \mod 7$: the closed (basic) Newton stratum is the union of all three non-$p$-ordinary \EO strata;
 \item if $p\equiv 1,6\mod 7$: the locally closed Newton stratum of codimension 1  is the union of two
 locally closed almost-$p$-ordinary \EO strata, and the closed (basic) Newton stratum agrees with the closed ($p$-special) \EO stratum.
 \end{enumerate}
   \end{enumerate}
\end{lemma}

By \cite[Theorem 1.6]{viehmann2013ekedahl}), for any PEL-type Shimura datum, the number of non-empty Ekedahl--Oort strata is independent of the prime $p$, while the number of non-empty Newton Polygon strata is not (it depends on the splitting behavior of $p$ in the associated reflex field). 
When $r=4$, the number of \EO strata of $\Sh(\D)_{\fpb}$ is equal to $2^N$, where $N=|\{\tau\in\tg\mid f(\tau)=1, f(\ts)=1\}|$.

Both the \EO types and Newton polygons corresponding to non-empty strata of depend on $p$. 
To illustrate this point, we list the Newton polygons associated with the Newton strata of $\Sh$;  we denote by $\mu$ the $\mu$-ordinary polygon and $\beta$ the basic polygon;  we write $(\lambda_1, \cdots \lambda_t)^s$ for the convex polygon with slopes $\lambda_1, \cdots \lambda_t$, each with multiplicity $s$.
\begin{enumerate}
    \item if $p\equiv 3,5$: $\mu=(\frac{1}{6}, \frac{5}{6})^6$,  and  $\beta=(\frac{1}{2}, \frac{1}{2})^6$;
    \item if $p\equiv 2,4$:   $\mu=(0,\frac{1}{3},\frac{2}{3},1)^3$, and  $\beta=(\frac{1}{3},\frac{2}{3})^6$;
    \item if $p\equiv 6\mod 7$:  $\mu=(0,1)^4 \oplus (\frac{1}{2}, \frac{1}{2})^2$,
       $\nu=(0,1)^2 \oplus (\frac{1}{2}, \frac{1}{2})^4$, and 
      $\beta= (\frac{1}{2}, \frac{1}{2})^6$;
    \item if $p\equiv 1\mod 7$:  $\mu=(0,1)^6$, $\nu= (0,1)^4\oplus (\frac{1}{2}, \frac{1}{2})^2$, and  $\beta= (0,1)^2\oplus (\frac{1}{2}, \frac{1}{2})^4$.   
    
\end{enumerate}


We deduce the following Lemma as a special case of \cite[Theorem 7.2]{lin2024abelian}  and Lemma \ref{lem:almost} combined.
\begin{lemma}
    \label{p6}\label{p1}
Assume $p \equiv 1, 6 \mod 7$ and $p>21$. For $\alpha \in \fpb\setminus\{0,1\}$,  let $C_\alpha$ be the corresponding smooth curve in $ \M_{\fpb}$. For $2\leq i\leq 3$, denote $\phi_i(t)=\phi_{\tau_i}(1,1)(t)\in \fpb [t]$ as defined in Proposition \ref{specialization of phi and psi}.  
Then
\begin{enumerate}
    \item if $\phi_2(\alpha)\cdot\phi_3(\alpha) \neq 0$, then $C_\alpha$ is $\D$-ordinary;
    \item if  either $\phi_2(\alpha)\neq 0$ and $\phi_3(\alpha)= 0$ or $\phi_2(\alpha)= 0$ and $\phi_3(\alpha)\neq 0$, then  $C_\alpha$ is almost $\D$-ordinary, each instance corresponding to a different \EO type; 
    \item if   $\phi_2(\alpha)=0$ and $\phi_3(\alpha)=0$, then $C_\alpha$ is $\D$-special.
\end{enumerate}
\end{lemma}
\begin{proof}
Recall the signature type of $\Sh$ is $(0,1,1,1,1,2)$; note that $g(\tau)=2$ for all $\tau\in \tg$. 
For $\tau\in \tg$ satisfying $f(\tau)=1$ (and $f(\tau^*)=1$), we denote by $h_\tau=h_{1,\tau}$ the $\tau$-partial Hasse--Witt invariant from Definition \ref{def:h_1}, and by $h_\OO=h_{0,\OO_\tau}$ the  $\OO$-partial Hasse--Witt invariant from Definition \ref{def:h_0}. 

Assume $p \equiv 6 \mod 7$. Then, $\tg$ splits into three Frobenius orbits: $\OO_i=\{\tau_i,\tau_i^*\}$, for $1\leq i\leq 3$. When $i=1$, $f(\tau_1)=0$, hence both the $\tau_1$-partial and $\OO_1$-partial Hasse--Witt invariants are trivial. When $2\leq i\leq 3$, $f(\tau_i)=1$ and $f(p\tau_i)=f(\ts_i)=1$, hence by definition $h_{\OO_i}=h_{\tau_i} h_{\ts_i}$ where $h_{\tau_i}=\phi_i$ and $h_{\ts_i}$ and $h_{\tau_i}$ have the roots in $\fp\setminus \{0,1\}$.


Assume $p \equiv 1 \mod 7$. Then, $\tg$ splits into six Frobenius orbits: $\OO_i=\{\tau_i\}$, for $1\leq i\leq 6$. When $i=1$ (resp. $i=6$) the $\tau$-partial Hasse--Witt invariant is trivial; since $h_{\ts_i}$ and $h_{\tau_i}$ have the roots in $\fp\setminus \{0,1\}$, it suffices to consider $2\leq i\leq 3$. By definition, for $2\leq i\leq 3$, $h_{\OO_{\tau_i}}=h_{\tau_i}=\phi_i$.
\end{proof}

\begin{proposition}\label{p equiv mo}
Assume $p \equiv 1,6 \mod 7$ and $p\neq 13$, then there exists $\alpha\in\fpb\setminus\{0,1\}$ satisfying $\phi_2(\alpha)=0$ and $\phi_3(\alpha) \neq 0$ (resp. $\phi_2(\alpha) \neq 0$ and $\phi_3(\alpha) =0$). 
\end{proposition}
\begin{proof}
Assume $p\equiv 6\pmod{7}$. 
By Proposition \ref{prop: separability of phit and psit}, the polynomials  $\phi_2(t)$ and $\phi_3(t)$ have only simple roots, except possibly $0,1$. From Notation \ref{define c C N}, $\phi_2=f(\frac{2p-5}{7},\frac{2p-5}{7},\frac{p-6}{7})$, hence by Lemma \ref{Lem: deg, vt and vt-1 of fabc},   $v_t(\phi_2)=v_{t-1}(\phi_2)=0$ and $\deg(\phi_2)=\frac{p-6}{7}$. Similarly,  $\phi_3=f(\frac{3p-4 }{7},\frac{3p-4 }{7},\frac{5p-2 }{7})$, and by Lemma \ref{Lem: deg, vt and vt-1 of fabc}, $v_t(\phi_3)=
\frac{2p+2}{7}$, $v_{t-1}(\phi_3)=0$ and $\deg(\phi_3)=\frac{3p-4}{7}$.
Set $\varphi_3={t^{-\frac{2p+2}{7}}}{\phi_3}$. Then, $\phi_2$ and $\varphi_3$ do not vanish at $0$ or $1$, 
have only simple roots and the same degree.
Hence, to conclude it suffices to show $\varphi_3 \neq \lambda \phi_2$, for some $\lambda\in \fpb$. 
By Lemma \ref{removing t and t-1},  $\varphi_3=f(\frac{3p-4}{7},\frac{3p-4}{7},\frac{p-6}{7})$. 
For $p\neq 13$,  $\min\{\frac{3p-4}{7},\frac{p-6}{7},2\frac{3p-4}{7}-\frac{p-6}{7}\}=\frac{p-6}{7}\geq 3$,  
and we deduce the statement from Lemma \ref{Lem f(abc) constant multiple}.

Assume $p\equiv 1\pmod{7}$. The argument is similar. From  $\phi_2=f(\frac{2p-2}{7},\frac{2p-2}{7},\frac{p-1}{7})$, we deduce $v_t(\phi_2)=v_{t-1}(\phi_2)=0$ and $\deg(\phi_2)=\frac{p-1}{7}$.
Similarly, from $\phi_3=f(\frac{3p-3}{7},\frac{3p-3}{7},\frac{5p-5}{7})$, we deduce $v_t(\phi_3)=\frac{2p-2}{7}$, $v_{t-1}(\phi_3)=0$ and $\deg(\phi_3)=\frac{3p-3}{7}$.
Hence,  $\phi_2$ and $\varphi_3=t^{-\frac{2p-2}{7}}{\phi_3}$ do not vanish at $0$ and $1$, have only simple roots and the same degree. Since $\varphi_3=f(\frac{3p-3}{7},\frac{3p-3}{7},\frac{p-1}{7})$, and $\min\{\frac{3p-3}{7},\frac{p-1}{7},2\frac{3p-3}{7}-\frac{p-1}{7}\}=\frac{p-1}{7}\geq 3$, 
we conclude by Lemma \ref{Lem f(abc) constant multiple}.
\end{proof}
\begin{lemma}\label{last}
  Assume $p \equiv 6 \mod 7$, then $\phi_2(-1)=\phi_3(-1)=0$.
\end{lemma}
\begin{proof}
 For $\alpha=-1$, the curve $C_\alpha$ has an extra automorphism of order $2$, and its Jacobian $J_\alpha$ has complex multiplication. By the Shimura-Tanayama formula, if $p \equiv 6 \mod 7$, then Newton polygon of $C_\alpha$ is $(\frac{1}{2}, \frac{1}{2})^6$, hence $C_\alpha$ is $\D$-special.
 
Alternatively, recall $\phi_2(t) =\sum_{i_2+i_3=\frac{p-6}{7}}\binom{\frac{2p-5}{7}}{i_2}\binom{\frac{2p-5}{7}}{i_3}t^{i_2}$. Hence, $\phi_2(t)$ is the coefficient of $x^{\frac{p-6}{7}}$ in the polynomial $(xt+1)^{\frac{2p-5}{7}}(x+1)^{\frac{2p-5}{7}}$, and $\phi_2(-1)$ is the coefficient of $x^{\frac{p-6}{7}}$ in $
(1-x^2)^{\frac{2p-5}{7}}$. 
Since $p$ is odd and $p \equiv 6 \pmod{7}$, $\frac{p-6}{7}$ is odd, hence $\phi_2(-1)=0$.
A similar argument show that $\phi_3(-1)=0$ since $\frac{5p-2}{7}$ is odd. 
\end{proof}

\begin{proof}[Proof of Proposition \ref{prop: special case}]
If $p\not\equiv 1, 6\pmod{7}$, the statement follows from Proposition \ref{prop: main theorem for r eq 4} for $p>21$, and by direct computations for $p<21$. 
If $p\equiv 1,6\pmod{7}$, $p\neq 13$:  the statement follows from Propositions \ref{p equiv mo} and Lemma \ref{last}.

For $p=13$, direct computations yield $\phi_2(t)=3(t+1)$ and $\phi_3(t)=5t^4(t+1)$. Thus $\phi_2(-1)=\phi_3(-1)=0$, and there is no $\alpha\in \fpb\setminus\{0,1\}$ which is a root of one of $\phi_2$, $\phi_3$ but not of the other.
\end{proof}


\begin{remark}\label{some p} Assume $p\equiv 1\pmod{7}$, and set $a=\frac{p-1}{7}$. 
By Lemma \ref{p1}, given $\alpha\in\fpb\setminus\{0,1\}$, $C_\alpha$ is $\D$-special  if and only if  $\phi_2(\alpha)=\phi_3(\alpha)=0$; by Proposition \ref{specialization of phi and psi}, $\phi_2=f(2a,2a;a)$ and $\phi_3= t^{2a} \varphi_3$, with $\varphi_3=f(3a,3a;a)$.
 By direct computations, $\alpha\in \fpb\setminus\{0,1\}$ satisfying $\gcd(\phi_2, \phi_3))(\alpha)=0$ exists for thirty two out of the first one hundred primes satisfying $p\equiv 1\bmod 7$ (in particular, it exists for $p=113$, but not for $p=29$).\footnote{
     If $p=29$ then $\phi_2(t)=12 (t^2 + t + 12)(t^2 + 17t + 17)$,  $\phi_3(t)=2t^8 (t^2 + 6t + 20) (t^2 + 9t + 16)$, and  $\gcd(\phi_2, \phi_3)=1$.
     If $p=113$ then $\phi_2(t)= 55 (t^2 + 42t + 1) (t^2 + 14t + 69) (t^2 + 87t + 95) (t^{10} + 40t^9 + 87t^8 + 61t^7 + 35t^6 + 91t^5 + 35t^4 + 61t^3 + 87t^2 + 40t + 1)$,  $\phi_3(t)=2 t^{32} (t^2 + 42t + 1)  (t^{14} + 84t^{13} + 34t^{12} + 99t^{11} + 15t^{10} + 102t^9 + 76t^8 + 12t^7 + 76t^6 + 102t^5 + 15t^4 + 99t^3 + 34t^2 + 84t + 1)$, and $\gcd(\phi_2,\phi_3)=t^2 + 42t + 1$.}
\end{remark}


By the same methods used to prove Proposition \ref{prop: special case}, we deduce the following results for other families of cyclic covers of degree $5$ and $7$.

    \begin{proposition} \label{prop: other fam}
    Let $(m,4,\au)$ be a cyclic monodromy datum,  and $p$ a rational prime, $p\nmid m$.
     We denote 
    by $\D$ the associated  PEL-type Shimura datum. Assume either $m=5$ and $p>15$,
    or $m=7$  and $p>21$. When $m=7$, we further assume 
    $\au$ is not equivalent to $(1,6,c,7-c)$, for some $1\leq c\le 6$. 
    
    Then there exists a smooth {almost-$\D$-ordinary} cover of $\Po$, defined over $\fpb$, with monodromy datum $\gamma$. 
    More precisely, each
{almost-$p$-ordinary} Ekedahl-Oort stratum in $\Sh(\D)_{\fpb}$ has non-empty intersection with the image of $\M(m,4,\au)$ under the Torelli map. 
Furthermore,  there is a unique Newton stratum  of codimension 1 in $\Sh(\D)_{\fpb}$, and it has non-empty intersection with the image of $\M(m,4,\au)$ under the Torelli map.


    \item  
   
\end{proposition}
\begin{proof}
Set $\M=\M(m,4,\au)$. For $m=5$: without loss of generality, we may assume that the signature type $\fb$ of $\gamma$ is either $(0,1,1,2)$ or $(1,1,1,1)$ (that is $\fb$ is equal to one of the two, modulo the action of the Galois group $\Z/5\Z^*$). In the first case,  the statement is proven in \cite[Proposition 5.8]{li2019newton}.
In the second case, the non-emptiness of the intersection of the image of $\M$ with the almost-$p$-ordinary \EO strata
 follows from
     the results in Appendix \ref{sec: separability}, by arguments analogous to those in the proof of Proposition \ref{p equiv mo}.  Furthermore, in this instance, there is unique Newton stratum of codimension 1 in $\Sh(\D)_{\fpb}$, and it is equal to the union of the {almost-$p$-ordinary} \EO strata.
     
     Similarly, for $m=7$: without loss of generality, we may assume that the signature type $\fb$ of $\gamma$ is either $(0,0,1,1,2,2)$ or $(0,1,1,1,1,2)$ or $(1,1,1,1,1,1)$. In the first case, the statement is proven in \cite[Proposition 5.8]{li2019newton}; in the second case, it follows form Proposition \ref{prop: special case}. The third case only occurs for inertia types $\au$ equivalent  to $(1,6,c,7-c)$, for some $1\leq c\leq 6$.
\end{proof}

\begin{remark}
   Consider the cyclic monodromy datum $(7,4,(1,6,c,7-c))$, for some $1\leq c\leq 6$, and a prime $p$, $p\neq 7$. Denote by $\D$ the associated PEL-type Shimura datum, and by $\fb$ its signature type. Then $\fb=(1,1,1,1,1,1)$; $\Sh(\D)_{\fpb}$
has three {almost-$p$-ordinary} \EO strata and a unique Newton stratum of codimension 1. Furthermore, if $p\not\equiv 3,5\bmod 7$, then the Newton stratum of codimension 1 contains all almost-$p$-ordinary \EO strata. 

For any prime $p>21$, there exists a smooth almost-$\D$-ordinary cover with the prescribed monodromy, and hence if $p\not \equiv 3,5 \bmod 7$ the intersection of the Newton stratum of codimension 1 with the image of $\M(7,4,(1,6,c,7-c))$ is non-empty.

In this instance, our methods enable us to identify one non-empty almost-$\D$-ordinary \EO stratum, its the choice depends on the congruence class $c\bmod 7$. 
\end{remark}

\bmhead{Acknowledgements}

We would like to thank Wanlin Li and Rachel Pries for many helpful discussions. Mantovan is partially supported by NSF grant DMS-22-00694.


\appendix


\section{The polynomials $f(a,b,c)$}\label{sec: separability}
In this section, we study the properties of the polynomials $f(a,b,c) \in \fp[t]$ from Definition \ref{def: fabc}, and prove Lemma  \ref{Lem: deg, vt and vt-1 of fabc}.
\begin{notation}
  For any integers $a,b,c\geq 0$ satisfying  $a,b<p$ and $c\leq a+b$, let $f(a,b,c)\in \fpb[t]$ be \[f(a,b,c)=\tsum_{i_2+i_3=c}\tbinom{a}{i_2}\tbinom{b}{i_3}t^{i_2}.\]  
\end{notation}

 \begin{lemma}[Lemma \ref{Lem: deg, vt and vt-1 of fabc} (3,4,5)] \label{lem:345}
\begin{enumerate}
     \item $\deg(f(a,b,c)) =\min\{a,c\}$;
 \item $ v_t(f(a,b,c)) =\max\{0,c-b\}$;
 \item $v_{t-1}(f(a,b,c)) = a+b-(p-1) $ if 
 $a+b-(p-1)\leq c\leq p-1<a+b$, and $0$ 
   otherwise.
 \end{enumerate}
 \end{lemma}
\begin{proof}
The first two statements 
follows from the definition; we compute $v_{t-1}(f(a,b,c))$. For $0\leq s\leq \min\{a,c\}$, let  $f(a,b,c)^{(s)}$ denote the $s$-{th} derivative of $f(a,b,c)$, then
$f(a,b,c)^{(s)}=\frac{a!}{(a-s)!} f(a-s,b,c-s)$ and  $f(a,b,c)^{(s)}(1)=\binom{a+b-s}{c-s}$. By assumptions, $a+b-s<2p$, and hence $v_p((a+b-s)!)\leq 1$.

If $a+b\leq p-1$, then $v_p((a+b)!)=0$; if $c\geq p$, then $v_p(c!)=1$; if $c\leq a+b-p$, then $p\leq a+b-c$  and $v_p((a+b-c)!)=1$. In each instance, $v_p(\binom{a+b}{c})=0$ and hence $f(a,b,c)(1)\neq 0$.

If $0<a+b-(p-1)\leq c\leq p-1$, then for any $0\leq s<a+b-(p-1)$, the inequalities $a+b-s>p-1\geq c-s\geq 0$ 
and $0\leq (a+b-s)-(c-s)\leq p-1$ imply $v_p(\binom{a+b-s}{c-s})=1$. For $s=a+b-(p-1)$, $v_p(\binom{a+b-s}{c-s})=0$.  Hence, $v_{t-1}(f(a,b,c)=a+b-(p-1)$.
\end{proof}
\begin{lemma}[Recurrence Relations]\label{Lem: fabc involution} 
\begin{enumerate}
    \item $f(a,b,c)+f(a,b,c-1)=f(a,b+1,c),$
    \item $f(a,b,c)+tf(a,b,c-1)=f(a+1,b,c),$
\item $f(a+1,b,c)-f(a,b+1,c)=(t-1)f(a,b,c-1),$
\item $\binom{a+b}{a} f(a,b,c) =\binom{a+b}{c} f(c,a+b-c,a)$.
\end{enumerate}
\end{lemma}
\begin{proof}
By definition, $f(a,b,c)$ is the coefficient of $\lambda^c$ in $(1+\lambda t)^a (1+\lambda)^b$. Hence, the first statement follows from identity $(1+\lambda t)^a (1+\lambda)^{b+1}= (1+\lambda t)^a (1+\lambda)^b +\lambda (1+\lambda t)^a (1+\lambda)^b$, the second statement from  $(1+\lambda t)^{a+1} (1+\lambda)^{b}= (1+\lambda t)^a (1+\lambda)^b +\lambda t(1+\lambda t)^a (1+\lambda)^b$, and 
the third statement is a consequence of the first two.

We verify the fourth statement: for $0\leq i\leq \max\{a,c\}$, the coefficients of $t^i$ in the two polynomials are
{
\setlength{\abovedisplayskip}{4pt}
\setlength{\belowdisplayskip}{4pt}
\setlength{\abovedisplayshortskip}{2pt}
\setlength{\belowdisplayshortskip}{2pt}
\setlength{\jot}{0pt}
\begin{align*}
\tbinom{a+b}{a}\tbinom{a}{i}\tbinom{b}{c-i}
&
=\tfrac{(a+b)!}{i!(a-i)! (c-i)! (b-c+i)!} 
=\tbinom{a+b}{c}\tbinom{c}{i}\tbinom{a+b-c}{a-i}. \qedhere
\end{align*}}
\end{proof}

\begin{lemma} [Extracting factors $t$ and $t-1$]\label{removing t and t-1}


\begin{enumerate}
    \item if $c\geq b$, then
$f(a,b,c)=t^{c-b}f(b,a,a+b-c);$
\item if 
$a+b-(p-1)\leq c\leq p-1\leq a+b$, then
{
\setlength{\abovedisplayskip}{4pt}
\setlength{\belowdisplayskip}{4pt}
\setlength{\abovedisplayshortskip}{2pt}
\setlength{\belowdisplayshortskip}{2pt}
\setlength{\jot}{0pt}
$$f(a,b,c)=(-1)^{a+c}(t-1)^{a+b-(p-1)} f(c-(a+b)+(p-1), (p-1)-c, (p-1)-b).$$}
\end{enumerate}
\end{lemma}
\begin{proof}
We verify the first statement: if $c\geq b$, by Lemma \ref{lem:345}, $f(a,b,c)$ is divisible by $t^{c-b}$; for $c-b\leq j\leq a$, the coefficients of $t^{j}$ in the two polynomials are 
{
\setlength{\abovedisplayskip}{4pt}
\setlength{\belowdisplayskip}{4pt}
\setlength{\abovedisplayshortskip}{2pt}
\setlength{\belowdisplayshortskip}{2pt}
\setlength{\jot}{0pt}
$$
\tbinom{b}{j-(c-b)}\tbinom{a}{(a+b-c)-(j-(c-b))}
= \tbinom{b}{c-j}\tbinom{a}{j}. 
$$}

We prove the second statement by induction on $s=a+b-(p-1)$, via  Lemma~\ref{Lem: fabc involution}.  
Recall the assumption $a+b-(p-1)\leq c\leq p-1\leq a+b$. If $s=0$, then
{
\setlength{\abovedisplayskip}{4pt}
\setlength{\belowdisplayskip}{4pt}
\setlength{\abovedisplayshortskip}{2pt}
\setlength{\belowdisplayshortskip}{2pt}
\setlength{\jot}{1pt}
\[\tbinom{p-1}{a}f(a,b,c)=\tbinom{p-1}{c}f(c,p-1-c,a)=\tbinom{p-1}{c}f(c,p-1-c,p-1-b),\]}
and $\binom{p-1}{k}\equiv (-1)^k\pmod{p}$.

If $s\geq 1$, then 
{
\setlength{\abovedisplayskip}{4pt}
\setlength{\belowdisplayskip}{4pt}
\setlength{\abovedisplayshortskip}{2pt}
\setlength{\belowdisplayshortskip}{2pt}
\setlength{\jot}{1pt}
\begin{align*}
f(a,b,c)&=f(a,b-1,c)+f(a,b-1,c-1)\\
&=(-1)^{a+c}(t-1)^{s-1} f(c-s+1,p-1-c,p-b) + (-1)^{a+c-1}(t-1)^{s-1} f(c-s,p-c,p-b)\\
&=(-1)^{a+c}(t-1)^{s-1}\big( f(c-s+1,p-1-c,p-b) - f(c-s,p-c,p-b) \big)\\
&= (-1)^{a+c}(t-1)^{s-1} (t-1) f(c-s,p-1-c,p-1-b). \qedhere
\end{align*}}
\end{proof}

\begin{lemma}\label{recurrance one} 

If  $a+b\leq p-1$ or $c\leq a+b-p+1$, then 
$\gcd(f(a,b,c), f(a,b,c-1))=1$.
\end{lemma}

\begin{proof}
We  prove the statement by induction on $c$.
If $c=1$, then $f(a,b,1)=b+at$ and $f(a,b,0)=1$. 

Assume $c\geq 2$. 
First, we establish the following recurrence relations:
{
\setlength{\abovedisplayskip}{4pt}
\setlength{\belowdisplayskip}{4pt}
\setlength{\abovedisplayshortskip}{2pt}
\setlength{\belowdisplayshortskip}{2pt}
\setlength{\jot}{1pt}
\begin{align}
    &c f(a,b,c)=(a-c+1)tf(a,b,c-1)+b f(a+1,b-1,c-1), \label{rel 1}\\
    &(a-c)f(a,b-1,c)=af(a-1,b,c)+(c-a-b)f(a,b-1,c-1), \label{rel 2}\\
    &(a+2-c)f(a+1,b-1,c-1)=(a+1)f(a,b,c-1)+(c-2-a-b)f(a+1,b-1,c-2).\label{rel 3}
\end{align}}
For \eqref{rel 1}: by definition, $tf(a,b,c-1)(t)=\sum_i \binom{a}{i}\binom{b}{c-1-i}t^{i+1}=\sum_i\binom{a}{i-1}\binom{b}{c-i}t^i$. 
Comparing the coefficients of $t^i$, 
{
\setlength{\abovedisplayskip}{4pt}
\setlength{\belowdisplayskip}{4pt}
\setlength{\abovedisplayshortskip}{2pt}
\setlength{\belowdisplayshortskip}{2pt}
\setlength{\jot}{0pt}
\begin{align*}
c\tbinom{a}{i}\tbinom{b}{c-i}&-(a-c+1)\tbinom{a}{i-1}\tbinom{b}{c-i}=
\tbinom{a}{i-1}\tbinom{b}{c-i}\big(\tfrac{c(a-i+1)}{i}-\tfrac{a-c+1}{1}\big)\\
&=\tbinom{a}{i-1}\tbinom{b}{c-i}\tfrac{(a+1)(c-i)}{i}=b \tbinom{a+1}{i}\tbinom{b-1}{c-1-i}.\end{align*}}

For \eqref{rel 2}: Similarly, comparing the coefficients of $t^i$, 
{
\setlength{\abovedisplayskip}{4pt}
\setlength{\belowdisplayskip}{4pt}
\setlength{\abovedisplayshortskip}{2pt}
\setlength{\belowdisplayshortskip}{2pt}
\setlength{\jot}{0pt}
\begin{align*}
(a-c)&\tbinom{a}{i}\tbinom{b-1}{c-i}-a\tbinom{a-1}{i}\tbinom{b}{c-i} =
\tbinom{a}{i}\tbinom{b-1}{c-i}\big((a-c)-a\tfrac{a-i}{a}\tfrac{b}{b-c+i} \big)\\
&=
\tbinom{a}{i}\tbinom{b-1}{c-i}\tfrac{(c-a-b)(c-i)}{(b-c+i)}= (c-a-b)\tbinom{a}{i}\tbinom{b-1}{c-1-i}. \end{align*}}
For \eqref{rel 3}: it follows from \eqref{rel 2} with  $a$ replaced by $a+1$ and $c$ by $c-1$.

By assumptions $1\leq b\leq p-1$, Equation \eqref{rel 1} implies 
{
\setlength{\abovedisplayskip}{4pt}
\setlength{\belowdisplayskip}{4pt}
\setlength{\abovedisplayshortskip}{2pt}
\setlength{\belowdisplayshortskip}{2pt}
\setlength{\jot}{1pt}
\[\gcd(f(a,b,c),f(a,b,c-1)) \mid \gcd(f(a,b,c-1), f(a+1,b-1,c-1)).\]}
From $a+b\leq p-1$, we deduce $ 3\leq a+2\leq a+b+2-c\leq a+b\leq p-1$;
    from $c\leq a+b-p$,  $ p+1\leq a+b+2-c\leq 2p-2$.
In both instance $(a+b+2-c)\not\equiv0\pmod{p}$.
Therefore, Equation \eqref{rel 3} implies 
{
\setlength{\abovedisplayskip}{4pt}
\setlength{\belowdisplayskip}{4pt}
\setlength{\abovedisplayshortskip}{2pt}
\setlength{\belowdisplayshortskip}{2pt}
\setlength{\jot}{1pt}
\[ \gcd(f(a,b,c-1), f(a+1,b-1,c-1)) \mid \gcd(f(a+1,b-1,c-1),f(a+1,b-1,c-2)) .\]}
Since $b\geq c\geq 2$,  the integers $a+1,b-1,c-1>0$ satisfy  $(a+1)+(b-1)=a+b\leq p-1$ and $c-1\leq b-1$. By inductive hypothesis, $\gcd(f(a+1,b-1,c-1),f(a+1,b-1,c-2))=1$, and we conclude  $\gcd(f(a,b,c),f(a,b,c-1))=1$.
\end{proof}


\begin{proposition}[Lemma \ref{Lem: deg, vt and vt-1 of fabc}(1)]\label{f(a,b,c) is separable} 
The polynomial $f(a,b,c)$ has only simple roots except possibly $0,1$.
\end{proposition}
\begin{proof}
Without loss of generality, we assume $a,b, c>0$ (as the remaining cases are trivial).
By definition, the derivative of $f(a,b,c)$ is
$f'(a,b,c)=af(a-1,b,c-1)(t)$; hence it suffices to show that $0,1$ are the only possible roots of the polynomial $\gcd(f(a,b,c),f(a-1,b,c-1))$.
By Lemma \ref{Lem: fabc involution}(2),  $f(a,b,c)=f(a-1,b,c)+t f(a-1,b,c-1)$; 
 we deduce that  
 {
\setlength{\abovedisplayskip}{4pt}
\setlength{\belowdisplayskip}{4pt}
\setlength{\abovedisplayshortskip}{2pt}
\setlength{\belowdisplayshortskip}{2pt}
\setlength{\jot}{1pt}
 \[\gcd(f(a,b,c),f(a-1,b,c-1))=\gcd(f(a-1,b,c-1),f(a-1,b,c)).\]}
By Lemma \ref{recurrance one}, if
$c\leq b$ and either $a+b\leq p$ or $c\leq a+b-p$, we conclude.


 By Lemma \ref{removing t and t-1},  there exist  integers $a_1,b_1,c_1\geq 0$ satisfying $a_1,b_1\leq p-1$ and $c_1\leq a_1+b_1$ such that
 $f(a,b,c)={t^{s_1}(t-1)^{s_2}}f(a_1,b_1,c_1)$, 
with $s_1=v_t(f(a,b,c))$ and $s_2=v_{t-1}(f(a,b,c))$. By construction,   $f(a_1,b_1,c_1)(0)\neq 0$ and $f(a_1,b_1,c_1)(1)\neq 0$, hence by Lemma \ref{Lem: fabc involution},  $c_1\leq b_1$, and  
either $a_1+b_1\leq p-1$ or $c_1\leq a_1+b_1-p$. In both instances, $f(a_1,b_1,c_1)$ has only simple roots, which suffices to conclude.
\end{proof}


The following lemma is the key ingredient in the proof of Propositions \ref{p equiv mo} and \ref{prop: other fam}.

\begin{lemma}\label{Lem f(abc) constant multiple} \label{lem: Nec cond for same roots}
For $1\leq i\leq 2$, let $a_i,b_i,c_i$ be non-negative integers 
satisfying $a_i,b_i<p$ and $c_i\leq a_i+b_i$. 
Assume $\min\{a_1,b_1,c_1,a_1+b_1-c_1\}\geq 3$.



If $f(a_1,b_1,c_2)$ and $f(a_2,b_2,c_2)$ share all roots except possibly $0,1$, then either $a_1+b_1=a_2+b_2$ or $a_1+b_1+a_2+b_2=2p-2$.
\end{lemma}
\begin{proof}

By  Lemma \ref{removing t and t-1}, it suffices to prove that
if $f(a_1,b_1,c_1)=\lambda f(a_2,b_2,c_2)$, for some $ \lambda\in\fp$, $\lambda\neq 0$, 
then $a_1+b_1=a_2+b_2$ and
$\{a_1,c_1,c_2-b_2-1\}\equiv\{a_2,c_2,c_1-b_1-1\}\bmod{p}$.

Let $v=v_t(f(a_1,b_1,c_1))=\max\{0,c_1-b_1\}$ and $d=\deg(f(a_1,b_1,c_1))=\min\{a_1,c_1\}$.   The equality $f(a_1,b_1,c_1)=\lambda f(a_2,b_2,c_2)$ implies, for any $v\leq i<d$,
{
\setlength{\abovedisplayskip}{4pt}
\setlength{\belowdisplayskip}{4pt}
\setlength{\abovedisplayshortskip}{2pt}
\setlength{\belowdisplayshortskip}{2pt}
\setlength{\jot}{1pt}
\[\tfrac{\binom{a_1}{i+1}\binom{b_1}{c_1-i-1}}{\binom{a_1}{i}\binom{b_1}{c_1-i}}=
\tfrac{(a_1-i)(c_1-i)}{(b_1-c_1+i+1)}=\tfrac{\binom{a_2}{i+1}\binom{b_2}{c_2-i-1}}{\binom{a_2}{i}\binom{b_2}{c_2-i}}
=\tfrac{(a_2-i)(c_2-i)}{(b_2-c_2+i+1)}.\]}

Let $g_1(x)=(a_1-x)(c_1-x)(b_2-c_2+x+1), g_2(x)=(a_2-x)(c_2-x)(b_1-c_1+x+1)\in \fp[x]$.  By definition, $g_1,g_2$ are monic, cubic polynomials satisfying $g_1(i)=g_2(i)$, for all $v\leq i<d$. By assumption $d-v\geq 3$,  hence $g_1=g_2$, that is  $\{a_1,c_1,c_2-b_2-1\}\equiv\{a_2,c_2,c_1-b_1-1\}\bmod{p}$. 
In particular, $a_1+b_1\equiv a_2+b_2\bmod{p}.$ By Lemma \ref{lem:345}, we deduce $a_1+b_1= a_2+b_2$.
\end{proof}

We conclude with a lemma on the binomial coefficients in characteristic $p$, which we use in the proof of Lemma \ref{aktphi}.
\begin{lemma}\label{choose p less}
For $p$ prime, and integer $c\geq d\geq p$,  we have
$\binom{c}{d}\equiv \binom{c}{d-p}\frac{1+\lfloor\frac{c-d}{p}\rfloor}{\lfloor\frac{d}{p}\rfloor} \mod{p}$.
\end{lemma}
\begin{proof}
Explicitly, $\tbinom{c}{d}=\tbinom{c}{d-p}\frac{(c-d+1) \cdots (c-d+p)}{(d-p+1) \cdots d}.$
Among the factors in the numerator, 
there is one multiple of $p$, that is $p(1+\lfloor\frac{c-d}{p}\rfloor)$, and the others  are congruent to $1,\dots,p-1$ mod $p$.
Similarly, among the factors in the denominator, there is one multiple of $p$, that is $p(\lfloor\frac{d}{p}\rfloor)$, and the  others are congruent to $1,\dots,p-1$ mod $p$.
\end{proof}

\section{Technical lemmas}\label{app: technical lemmas} 

This appendix contains explicit polynomial computations, and specifically the proof of Lemma~\ref{lem: entry of phi and psi},  used in the proof Proposition \ref{main proposition both} to bound the order of vanishing of the $\D$-ordinary $\OO$-partial Hasse--Witt invariant along the boundary divisor $x_r-x_{r-1}=0$.

First, we introduce notations for two classes of auxiliary polynomials $f$ and $q$ in $\fpb[x_1, \dots x_r]$, which we use to express the 
entries of the maps $\phi_\tau$ and $\psi_\tau$ (see Remark \ref{phipsiasf}). 
In Lemma \ref{derivative and zero}, we establish recursive formulas for the derivatives of $f$ and $q$, with respect to the variable $x=x_r-x_{r-1}$ that we use to prove Lemma~\ref{lem: entry of phi and psi}.

\begin{notation}\label{define c C N}
For $ b_1, b_2, \dots, b_r,N, a, j' \in \mathbb{N}_{\geq 0}$, with $1 \leq j' \leq r-1$, $0 \leq a \leq r-1$ and $N< \sum_{k=1}^r b_k$, we define
{
\setlength{\abovedisplayskip}{4pt}
\setlength{\belowdisplayskip}{4pt}
\setlength{\abovedisplayshortskip}{2pt}
\setlength{\belowdisplayshortskip}{2pt}
\setlength{\jot}{1pt}
\begin{align*}
    c\bun&=\min\{1\leq c\leq r : \tsum_{i=1}^c b_i >N\}, \\
    C\bun&=N-\tsum_{i=1}^{c-1} b_i, \text{ here } c=c\bun,\\
    X\bun &=x_1^{b_1} \dots x_{c-1}^{b_{c-1}}x_c^C, \text{ here } c=c\bun, C=C\bun,\\
    f\bun &= \tsum_{i_1+ \dots +i_r=N} \tbinom{b_1}{i_1} \dots \tbinom{b_r}{i_r} x_1^{i_1} \dots x_r^{i_r}, \\
    q(a,k,r) &= -\tsum_{\substack{1 \leq i_1 <i_2< \dots < i_a \leq r\\ \text{each }i_t\neq k}} x_{i_1} \dots x_{i_a} \ ; \ q(0,k,r)=-1,\\
    \psi\bunj &=\tsum_{k=1}^r b_k f(b_1, \dots, b_{k}-1, \dots, b_r; N)q(r-j',k,r),\\
    s\bu &=\tsum_{k=1}^r b_i. 
\end{align*}}
\end{notation}
Notice that by definition, $1\leq c\bun \leq r$, $0 \leq C\bun < b_c$. 
\begin{remark}\label{phipsiasf}
For $G$ cyclic of size $m$, and $\tau \in \tg$,
if we set $b_{k} = \ptk$, for all $1\leq k\leq r$, and $s_\tau=s(b_1, \cdots, b_r)=\sum_{k=1}^r \ptk$, we have
{
\setlength{\abovedisplayskip}{4pt}
\setlength{\belowdisplayskip}{4pt}
\setlength{\abovedisplayshortskip}{2pt}
\setlength{\belowdisplayshortskip}{2pt}
\setlength{\jot}{1pt}
\begin{align*}
\phit(j',j)&=(-1)^{\st-pj+j'}f(b_1, \dots, b_r; \st-pj+j'), \text{ and }\\
\psit(j',j) &=(-1)^{\st-pj+r-j'}\psi(b_1,\dots,b_r;\st-pj;j')\\
&=(-1)^{\st-pj+r-j'}\tsum_{k=1}^r b_k f(b_1, \dots, b_{k}-1, \dots, b_r; \st-pj)q(r-j',k,r).
\end{align*}}
\end{remark}

In the following, given $g=g(x_1,\dots,x_r)\in \fpb[x_1, \dots, x_{r}]$, we  define the polynomial $g(x)$ by setting $x_r=x_{r-1}+x$, that is $g(x)=g(x_1,\dots,x_{r-1}, x_{r-1}+x)$. We denote by $g(0)$ the evaluation of $g(x)$ at $x=0$, hence $g(0)\in\fpb[x_1, \dots, x_{r-1}]$ .


  

\begin{lemma}\label{derivative and zero}
Consider the polynomials $f$ and $q$ as in Notation \ref{define c C N}. Then, we have
\begin{enumerate}
    \item $\frac{\partial^{t}}{\partial x^{t}}f\bun(x)=t!\binom{b_r}{t}f(b_1, \dots, b_{r-1},b_r-t; N-t)(x)$ for all $t\geq 0$;
    \item $f\bun(0)=f(b_1, \dots,b_{r-2}, b_{r-1}+b_r; N)$;
    \item $\frac{\partial}{\partial x}q(a, k, r)(x)=q(a-1,k,r-1)$ if $ k \leq r-1$ and $a \geq 1$, and $\frac{\partial}{\partial x}q(a, k, r)(x)=0$ if either $k=r$ or $a=0$.
    
    \item $\frac{\partial^2}{\partial x^2}q(a, k, r)(x) = 0$;    
    \item for every non-zero monomial $M$ of $q(a,k,r)(0)$: $ v_{x_s}(M) \leq 2$, for all $1\leq s\leq r-1$.   
\end{enumerate}
\end{lemma}
\begin{proof}
\begin{enumerate}
\item We induct on $t$, starting from $t=0$. The base case is evident. Assume the statement holds for some $t\geq 0$. Then
{
\setlength{\abovedisplayskip}{4pt}
\setlength{\belowdisplayskip}{4pt}
\setlength{\abovedisplayshortskip}{2pt}
\setlength{\belowdisplayshortskip}{2pt}
\setlength{\jot}{1pt}
\begin{align*}
&\tfrac{\partial^{t+1}}{\partial x^{t+1}}f\bun(x)
=t!\tbinom{b_r}{t} \tfrac{\partial}{\partial x} f(b_1, \dots, b_{r-1},b_r-t; N-t)(x)\\
&=t!\tbinom{b_r}{t} \frac{\partial}{\partial x} 
\tsum_{i_1+\dots+i_r=N-t} 
\tbinom{b_1}{i_1}\dots \tbinom{b_{r-1}}{i_{r-1}} \tbinom{b_r-t}{i_r} x_1^{i_1}\dots x_{r-1}^{i_{r-1}} (x+x_{r-1})^{i_r}\\
&=t!\tbinom{b_r}{t} 
\tsum_{i_1+\dots+i_r=N-t} 
\tbinom{b_1}{i_1}\dots \tbinom{b_{r-1}}{i_{r-1}} 
\tfrac{b_r-t}{i_r}\tbinom{b_r-t-1}{i_r-1} 
x_1^{i_1}\dots x_{r-1}^{i_{r-1}} i_r (x+x_{r-1})^{i_r-1}\\
&=(t+1)!\tbinom{b_r}{t+1} f(b_1, \dots, b_{r-1},b_r-t-1; N-t-1)(x).
\end{align*}}

\item  By definition, 
{
\setlength{\abovedisplayskip}{4pt}
\setlength{\belowdisplayskip}{4pt}
\setlength{\abovedisplayshortskip}{2pt}
\setlength{\belowdisplayshortskip}{2pt}
\setlength{\jot}{1pt}
\begin{align*}
f\bun(0)
&=\tsum_{i_1+\dots+i_r=N} 
\tbinom{b_1}{i_1}\dots
\tbinom{b_{r-1}}{i_{r-1}} 
\tbinom{b_r}{i_r} x_1^{i_1}
\dots x_{r-1}^{i_{r-1}} (0+x_{r-1})^{i_r}\\
&=\tsum_{i_1+\dots+i_{r-2}+i=N} 
\tbinom{b_1}{i_1}\dots 
\tbinom{b_{r-2}}{i_{r-2}}  
\big(\tsum_{i_{r-1}+i_r=i}
\tbinom{b_{r-1}}{i_{r-1}} \tbinom{b_r}{i_r}\big)
x_1^{i_1}\dots x_{r-2}^{i_{r-2}}x_{r-1}^{i}  \\
&=\tsum_{i_1+\dots+i_{r-2}+i=N} 
\tbinom{b_1}{i_1}\dots 
\tbinom{b_{r-2}}{i_{r-2}} 
\tbinom{b_{r-1}+b_r}{i} x_1^{i_1}
\dots x_{r-2}^{i_{r-2}}x_{r-1}^{i}.  
\end{align*}}

\item By definition, if $k\leq r-1$ and $a \geq 1$, then $q(a, k, r)$ is linear in $x$, and the coefficient of $x$ is $q(a-1, k, r-1)$; if $k=r$ or $a=0$, then $q(a,r,r)\in \fp[x_1,\dots,x_{r-1}]$,  and hence $\frac{\partial}{\partial x}q(a,r,r)(x)=0$.

\item By part (3),  $\frac{\partial}{\partial x}q(a,k,r)(x)\in \fp[x_1,\dots,x_{r-1}]$, hence  $\frac{\partial^2}{\partial x^2}q(a,k,r)(x)=0$. 

\item This follows from the definition of $q(a,k,r)$. \qedhere
\end{enumerate}
\end{proof}

\begin{lemma}\label{lem: entry of phi and psi}
Let $(m,r,\au)$ be a cyclic monodromy datum and $p$ be a prime such that $p>m(r-2)$. 
Assume there exists $\tn \in \tg$ satisfying $\ker(\tn)=0$ and $f(\tns)=1$, and let $0 \leq i <l=l(\OO_{\tn})$.  
For any $k$, $1 \leq k \leq r$, define
{
\setlength{\abovedisplayskip}{4pt}
\setlength{\belowdisplayskip}{4pt}
\setlength{\abovedisplayshortskip}{2pt}
\setlength{\belowdisplayshortskip}{2pt}
\setlength{\jot}{1pt}
\[
b_k =
\begin{cases}
    {\lfloor{p\langle\tfrac{p^i\tau_0 a(k)}{m}\rangle}\rfloor} 
    & \text{ if $f(p^i\tn^*)\geq 1$} \\
    {\lfloor{p\langle\tfrac{-p^i\tn a(k)}{m}\rangle}\rfloor} 
    & \text{if $f(p^i\tn^*)=0$}. \\
\end{cases}
\]}
Let $J \in \mathfrak{J}$, and  $R_{J,i}\in \fp[x_1, \dots, x_r]$ be as in Equation \eqref{RJi}, and define $R_{J,i}(x) \in \fp[x_1,\dots,x_{r-1}][x]$  by setting $x_r=x_{r-1}+x$. 

\begin{enumerate}
    \item Assume either $ f(p^i\tn^*)\geq 1, f(p^{i+1}\tn^*)\geq 1$ or $ f(p^i\tn^*)=0, f(p^{i+1}\tn^*)=0$.  Then, for $t = \max\{0, b_r+b_{r-1}-(p-1)\}$, we have
    \begin{enumerate}
        \item $(\frac{\partial^{t}}{\partial x^{t}}R_{J,i}) (0)\in \fp[x_1,\dots,x_{r-1}]$ is not identically 0;
\item for every non-zero monomial $M$ in $(\frac{\partial^{t}}{\partial x^{t}}R_{J,i})(0)$: 
 $v_{x_k}(M)<p$ for all $1 \leq k \leq r-1$.
    \end{enumerate}
    \item Assume either $ f(p^i\tn^*)\geq 1, f(p^{i+1}\tn^*)= 0$ or $ f(p^i\tn^*)=0, f(p^{i+1}\tn^*) \geq 1$. Then, for $t = \max\{0, b_r+b_{r-1}-(p-3)\}$, we have  
    \begin{enumerate}
    \item $(\frac{\partial^{t}}{\partial x^{t}}R_{J,i}) (0)\in \fp[x_1,\dots,x_{r-1}]$ is not  identically 0 if $J(i)=J(i+1)=1$;
\item for every non-zero monomial $M$ of $(\frac{\partial^{t}}{\partial x^{t}}R_{J,i})(0)$:  $v_{x_k}(M)<p$ for all $1 \leq k \leq r-1$.
\end{enumerate}
\end{enumerate}
\end{lemma}
Note that  $R_{J,i}=\phi_{\tau}(J(i+1), J(i))$ in case (1), and $R_{J,i}=\psi_{\tau}(J(i+1), J(i))$ in case (2). 
\begin{proof}  Write $\tau=p^i \tn$ if $f(p^i\tn^*)\geq 1$, 
and $\tau=p^i \tn^*$ if $f(p^i\tn^*)=0$. By definition, for all $1\leq k\leq r$, $b_k=\ptk$. 
Set $s=\sum_{k=1}^r b_k$ (it depends on $i$). 

\begin{enumerate}
\item Suppose $ f(p^i\tn^*)\geq 1, f(p^{i+1}\tn^*)\geq 1$ or $ f(p^i\tn^*)=0, f(p^{i+1}\tn^*)=0$. 




Denote 
$N(J,i)=s-pJ(i)+J(i+1)$.
Then, 
{
\setlength{\abovedisplayskip}{4pt}
\setlength{\belowdisplayskip}{4pt}
\setlength{\abovedisplayshortskip}{2pt}
\setlength{\belowdisplayshortskip}{2pt}
\setlength{\jot}{1pt}
\[
(\tfrac{\partial^{t}}{\partial x^{t}}R_{J,i}) (0)
= t!\tbinom{b_r}{t} 
\tsum_{i_1+ \dots+ i_{r-1}=N(J,i)-t}
\tbinom{b_1}{i_1} \dots 
\tbinom{b_{r-1}+b_r-t}{i_{r-1}}x_1^{i_1} 
\dots x_{r-1}^{i_{r-1}}.
\]}

\begin{enumerate}
    \item We show that $(\frac{\partial^{t}}{\partial x^{t}}R_{J,i}) (0)$ is not identically $0$.  Since $b_{r-1},b_r\leq p-1$, we have $0\leq t\leq b_r\leq p-1$, and thus $t!\binom{b_r}{t}\neq 0$.
Hence,  
it suffices to show that there exists a choice of integers $i_1, \dots i_{r-1}$ satisfying $0\leq i_k\leq b_k$ for $1\leq k \leq r-2$,  $0\leq i_{r-1}\leq b_{r-1}+b_r-t$ and $\sum_{k=1}^{r-1} i_k=N(J,i)-t$. Since the sum of the upper bounds for $i_k$ is $\sum_{k=1}^{r-2}b_i+(b_{r-1}+b_r)-t=s-t$ and the sum of lower bounds is $0$, 
such a choice for $i_k's$ exists if and only if $0 \leq N(J,i)-t\leq s-t$. 

We prove $N(J,i) \geq t$. By assumption, $f(\ts) \geq 1$ and $f(p\ts) \geq 1$; by definition,
$1 \leq J(i) \leq d(\tau)=f(\ts)$ and $1 \leq J(i+1) \leq d(p\tau)=f(p\ts)$. Hence, 
{
\setlength{\abovedisplayskip}{4pt}
\setlength{\belowdisplayskip}{4pt}
\setlength{\abovedisplayshortskip}{2pt}
\setlength{\belowdisplayshortskip}{2pt}
\setlength{\jot}{0pt}
\[N(J,i)=p(f(\ts)-J(i)+1)-(f(p\ts)-J(i+1)+1) \geq p-f(p\ts) \geq p-(r-2).\] 
On the other hand, by assumption $p \geq m(r-2)$, and hence
\[
t=b_r+b_{r-1}-(p-1) 
\leq b_r
=\pbtar \leq p\tfrac{m-1}{m}=p-\tfrac{p}{m} \leq p-(r-2).\] }

We prove $N(J,i) \leq s$. It suffices to show $pJ(i)-J(i+1) \geq 0$. Since $J(i) \geq 1$ and $J(i+1) \leq f(p\ts) \leq r-2$, it follows from  the assumption $p \geq m(r-2)$.

\item We show that $v_{x_k}(M)<p$ for all non-zero monomial $M$ in $(\frac{\partial^{t}}{\partial x^{t}}R_{J,i})(0)$ and for all $1 \leq k \leq r-1$. For $1 \leq k \leq r-2$, the inequality follows from $i_k \leq b_k\leq p-1$; for $k=r-1$, the inequality follows from $i_{r-1} \leq b_{r-1}+b_r-t \leq p-1$.

\end{enumerate}

\item Suppose either $ f(p^i\tn^*)\geq 1, f(p^{i+1}\tn^*)=0$ or $ f(p^i\tn^*)=0, f(p^{i+1}\tn^*)\geq 1$.

Denote $N=N(J,i)=s-pJ(i)$. For convenience,  write $j=J(i)$ and $j'=J(i+1)$.
Then, with notations as in Notation \ref{define c C N}, by Remark \ref{phipsiasf}, we have
{
\setlength{\abovedisplayskip}{4pt}
\setlength{\belowdisplayskip}{4pt}
\setlength{\abovedisplayshortskip}{2pt}
\setlength{\belowdisplayshortskip}{2pt}
\setlength{\jot}{1pt}
\begin{align*}
\psit(j',j)(x) &= (-1)^{s-pj+r-j'}\psi\bunj(x)\\
&=(-1)^{s-pj+r-j'}\tsum_{k=1}^r b_k f(b_1, \dots, b_{k}-1, \dots, b_r; N)(x)q(r-j',k,r)(x).
\end{align*}}
We argue separately the cases $t=0$ and  $t>0$.

Suppose  $t=0$, that is $b_r+b_{r-1}\leq p-3$. Then
{
\setlength{\abovedisplayskip}{4pt}
\setlength{\belowdisplayskip}{4pt}
\setlength{\abovedisplayshortskip}{2pt}
\setlength{\belowdisplayshortskip}{2pt}
\setlength{\jot}{0pt}
\begin{align*}
&(-1)^{s-pj+r-j'} (\tfrac{\partial^{t}}{\partial x^{t}}R_{J,i}) (0)
=(-1)^{s-pj+r-j'}\psit(j',j)(0)\\
&=\tsum_{k=1}^r b_k f(b_1, \dots, b_{k}-1, \dots, b_r; N)(0)\cdot q(r-j',k,r)(0)=S_1+S_2,
\end{align*}}
where $S_2=(b_{r-1}+b_r) f(b_1,\dots, b_{r-2}, b_{r-1}+b_r-1;N) \cdot q(r-j',r,r)
$, and 
{
\setlength{\abovedisplayskip}{4pt}
\setlength{\belowdisplayskip}{4pt}
\setlength{\abovedisplayshortskip}{2pt}
\setlength{\belowdisplayshortskip}{2pt}
\setlength{\jot}{0pt}
\begin{align*}
S_1=\tsum_{k=1}^{r-2} b_k f(b_1, &\dots, b_{k}-1, \dots, b_{r-2},b_{r-1}+ b_r; N) \\
&\cdot \big( q(r-j',k,r-1) + x_{r-1} q(r-j'-1,k,r-1) \big).
\end{align*}}
Note that if $j=j'=1$, then, for any $1\leq k\leq r-2$, we have 
$q(r-j',k,r-1)=0$, $x_{r-1} q(r-j'-1,k,r-1)= -x_{r-1} \frac{x_1\dots x_{r-1}}{x_k},$
and $q(r-j',r,r)=-x_1\dots x_{r-1}.$

In particular, $v_{x_{r-1}} ( q(r-j',k,r-1) + x_{r-1} q(r-j'-1,k,r-1) )=2$.



\begin{enumerate}
\item We show that $(\frac{\partial^{t}}{\partial x^{t}}R_{J,i}) (0)$ is not identically zero when $j=j'=1$. In this case, $N=s -p \leq \tsum_{k=1}^r b_k -(b_{r-1}+b_r+3)<\tsum_{k=1}^{r-2}b_k$,
and thus  $c=c(b_1,\dots,b_{r-2},b_{r-1}+b_{r}-1;N)\leq r-2$.
Consider the monomial $Y=x_1\dots x_{r-1} X(b_1,\dots,b_{r-2},b_{r-1}+b_{r}-1;N)$.
Note that $v_{x_{r-1}}(Y)=1$. Since for all $1\leq k\leq r-2$,
{
\setlength{\abovedisplayskip}{4pt}
\setlength{\belowdisplayskip}{4pt}
\setlength{\abovedisplayshortskip}{2pt}
\setlength{\belowdisplayshortskip}{2pt}
\setlength{\jot}{0pt}
\[x_{r-1}^2\mid ( q(r-j',k,r-1) + x_{r-1} q(r-j'-1,k,r-1) )\]}
by the above computation, we deduce that 
the coefficient of $Y$ in $(\frac{\partial^{t}}{\partial x^{t}}R_{J,i}) (0)$ is 
$(-1)^{s-pj+r-j'}(b_{r-1}+b_{r})\tbinom{b_c}{C}$, where $C=C(b_1,\dots,b_{r-2},b_{r-1}+b_{r}-1;N)$. In particular, it is non-zero.

\item We show all nonzero monomials $M$ in $(\frac{\partial^{t}}{\partial x^{t}}R_{J,i}) (0)=S_1+S_2$ has $v_{x_k}(M)<p$ for all $1 \leq k \leq r-1$.  By Notation \ref{define c C N},
{
\setlength{\abovedisplayskip}{4pt}
\setlength{\belowdisplayskip}{4pt}
\setlength{\abovedisplayshortskip}{2pt}
\setlength{\belowdisplayshortskip}{2pt}
\setlength{\jot}{0pt}
\[
\deg_{x_k}f(b_1,\dots,b_k-1,\dots,b_{r-2},b_{r-1}+b_r;N)
\leq
\begin{cases}
b_k\leq p-3 & \text{if }1\leq k\leq r-2,\\
b_r+b_{r-1}\leq p-3 & \text{if }k=r-1.
\end{cases}
\]}
Hence, if $M$ is a non-zero monomial in 
$S_1$, then $v_{x_k}(M) \leq p-1$, for all $1 \leq k \leq r-1$.
By a similar argument, if $M $ is a non-zero monomial in $S_2$, then $v_{x_k}(M) \leq p-1$, for all $1 \leq k \leq r-1$. 
\end{enumerate}
This finishes the proof of Lemma \ref{lem: entry of phi and psi}(2) in the case $t=0$.

 Suppose $t>0$, that is $b_r+b_{r-1}> p-3$. Then, by Lemma \ref{derivative and zero},
 {
\setlength{\abovedisplayskip}{4pt}
\setlength{\belowdisplayskip}{4pt}
\setlength{\abovedisplayshortskip}{2pt}
\setlength{\belowdisplayshortskip}{2pt}
\setlength{\jot}{1pt}
\begin{align*}
  &(-1)^{s-pj+r-j'}
  \tfrac{\partial^{t}}{\partial x^{t}}R_{J,i} (x)
  =(-1)^{s-pj+r-j'}\tfrac{\partial^{t}}{\partial x^{t}}\psi\bunj(x)\\
  &= \tsum_{k=1}^r b_k \tfrac{\partial^{t}}{\partial x^{t}}f(b_1, \dots, b_{k}-1, \dots, b_r; N)(x)\cdot q(r-j',k,r)(x) \\
  &+t\tsum_{k=1}^r b_k \tfrac{\partial^{t-1}}{\partial x^{t-1}}f(b_1, \dots, b_{k}-1, \dots, b_r; N)(x)\cdot \tfrac{\partial}{\partial x}q(r-j',k,r)(x) \\
  &= t ! \tbinom{b_r}{t}\tsum_{k=1}^r b_k f(b_1, \dots, b_{k}-1, \dots, b_r-t; N-t)(x)\cdot q(r-j',k,r)(x)
  \\
  &+t!\tbinom{b_r}{t-1}\tsum_{k=1}^{r-1} b_k f(b_1, \dots, b_{k}-1, \dots, b_r-t+1; N-t+1)(x)\cdot q(r-j'-1,k,r-1).
\end{align*}}
Setting $x=0$, we deduce, up to a sign, 
$
 \frac{\partial^{t}}{\partial x^{t}}R_{J,i} (0)
=S_1+S_2+S_3+S_4 $ where
{
\setlength{\abovedisplayskip}{4pt}
\setlength{\belowdisplayskip}{4pt}
\setlength{\abovedisplayshortskip}{2pt}
\setlength{\belowdisplayshortskip}{2pt}
\setlength{\jot}{0pt}
\begin{align*}\label{psi sum}
  &S_1= t ! \tbinom{b_r}{t}\tsum_{k=1}^{r-2} b_k f(b_1, \dots, b_{k}-1, \dots, b_{r-1}+b_r-t; N-t)\cdot q(r-j',k,r),\\
  &S_2=t!\tbinom{b_r}{t}(b_{r-1}+b_r)f(b_1, \dots, b_{r-2},b_{r-1}+b_r-t-1;N-t)\cdot q(r-j',r,r) ,\\
  &S_3=t!\tbinom{b_r}{t-1}\tsum_{k=1}^{r-2} b_k f(b_1, \dots, b_{k}-1, \dots, b_{r-2},p-2 ; N-t+1)\cdot q(r-j'-1,k,r-1) ,\\
  &S_4= t!\tbinom{b_r}{t-1}b_{r-1}f(b_1, \dots, b_{r-2},p-3; N-t+1)\cdot q(r-j'-1,r-1,r-1).
\end{align*}}
\begin{enumerate}
\item We show that $(\frac{\partial^{t}}{\partial x^{t}}R_{J,i}) (0)$ is not identically zero when $j=j'=1$. In this case, $N-t+1=s-p-b_{r-1}-b_r+p-2=\sum_{k=1}^{r-2}b_k-2.$
Consider the monomial $Y=x_1 \dots x_{r-2}X(b_1, \dots, b_{r-2}; N-t+1)$. 

Note that $v_{x_{r-1}}(Y)=0$. By definition, for any $1\leq k\leq r-2$, we have $x_{r-1}^2\mid q(r-1,k,r)(0)$, $x_{r-1}\mid q(r-2,k,r-1)$, and $x_{r-1}\mid q(r-1,r,r)$; also, $q(r-2,r-1,r-1)=-x_{1}\dots x_{r-2}$. By the above computation, we deduce that the coefficient of $Y$ $(\frac{\partial^{t}}{\partial x^{t}}R_{J,i}) (0)$ is 
$t!\tbinom{b_r}{t-1}b_{r-1} \tbinom{b_{c}}{C},$
where $C=C(b_1,\dots,b_{r-2};N-t+1)$; in particular, it is non-zero.
\item  We show all nonzero monomials $M$ in $(\frac{\partial^{t}}{\partial x^{t}}R_{J,i}) (0)$ has $v_{x_k}(M)<p$ for all $1 \leq k \leq r-1$.  Let $M$ be a non-zero monomial  in $ \frac{\partial^{t}}{\partial x^{t}}R_{J,i} (0)$.
By the above formula, for $1 \leq k \leq r-1$, $v_{x_k}(M)$ is at most equal to the sum of the degrees in $x_k$ of some polynomials $f, q$ as above. From Notation \ref{define c C N}, if  $1 \leq k \leq r-2$, then $\deg_{x_k}(f) \leq b_k \leq p-3$,  $\deg_{x_k}(q) \leq 1$, and hence $v_{x_k}(M) \leq p-1$.
If $k=r-1$, then $\deg_{x_{r-1}}(f)\leq p-2$,  $\deg_{x_{r-1}}(q) \leq 1$, and hence $v_{x_{r-1}}(M) \leq p-1$.
(More precisely, $M$ is a monomial in $S_i$, $1\leq i\leq 4$, then $\deg_{x_{r-1}}(f)$ is at most
 $b_{r-1}+b_r-t \leq p-3 $ if $i=1$, at most
 $b_{r-1}+b_r-t-1 \leq p-4 $ if $i=2$, at most  $p-2$ if $ i=3$, and at most $ p-3 $ if $i=4$.)
\end{enumerate}
This finishes the proof of Lemma \ref{lem: entry of phi and psi}(2) in the case $t>0$. \qedhere
\end{enumerate}
\end{proof}

\bibliography{sn-bibliography}

\end{document}